\documentclass[11pt]{article}
\usepackage{epsfig}

\def\be{\begin{equation}}
\def\ee{\end{equation}}
\def\bea{\begin{eqnarray}}
\def\eea{\end{eqnarray}}
\def\bes{\begin{eqnarray*}}
\def\ees{\end{eqnarray*}}

\def\nn{\nonumber}
\def\lb{\label}
\def\bs{\setminus}

\def\vs{{\varsigma}}

\def\R{{\bf R}}
\def\C{{\bf C}}
\def\Z{{\bf Z}}
\def\K{{\bf K}}
\def\N{{\bf N}}
\def\U{{\bf U}}

\def\Q{{\bf Q}}

\def\aa{{\alpha}}
\def\bb{{\beta}}
\def\ga{{\gamma}}

\def\ka{{\kappa}}
\def\th{{\theta}}
\def\om{{\omega}}
\def\Om{{\Omega}}
\def\ep{{\epsilon}}
\def\lm{{\lambda}}
\def\Lm{{\Lambda}}

\def\sg{{\sigma}}
\def\dm{{\diamond}}

\def\vf{{\varphi}}

\def\K{{\cal K}}
\def\P{{\cal P}}

\def\M{{\cal M}}

\def\rank{{\rm rank}}
\def\Sp{{\rm Sp}}
\def\CG{{\rm CG}}
\def\mod{{\rm mod}}
\def\ol{\overline}
\def\td#1{\tilde{#1}}

\def\hb{\vrule height0.18cm width0.14cm $\,$}

\title{The index quasi-periodicity and multiplicity of closed geodesics}

\author{Huagui Duan$^{1}$\thanks{Partially supported by NNSF Grant 10801079,
RFDP Grant 200800551002, LPMC of MOE of China and Nankai University.
E-mail: duanhg@nankai.edu.cn} \quad and \quad Yiming
Long$^{2}$\thanks{Partially supported by the 973 Program of MOST,
NNSF, MCME, RFDP, LPMC of MOE of China, and Nankai University.
E-mail: longym@nankai.edu.cn}\\\\
$^{1}$ School of Mathematics\\
$^{2}$ Chern Institute of Mathematics and LPMC\\
Nankai University, Tianjin 300071\\ The People's Republic of China\\
}

\begin{document}
\date{May 16, 2010}
\maketitle

\begin{abstract}
{\it In this paper, we prove the existence of at least two distinct
closed geodesics on every compact simply connected irreversible or reversible
Finsler (including Riemannian) manifold of dimension not less than $2$. }
\end{abstract}

{\bf Key words}: Closed geodesics, quasi-periodicity, multiplicity,
compact simply connected manifolds.

{\bf 2000 Mathematics Subject Classification}: 53C22, 58E05, 58E10.

\renewcommand{\theequation}{\thesection.\arabic{equation}}
\renewcommand{\thefigure}{\thesection.\arabic{figure}}

\setcounter{equation}{0}
\section{Introduction and main results}

The closed geodesic problem is a traditional and active topic in dynamical
systems and differential geometry for more than one hundred years.
Studies of closed geodesics can be traced back to J. Jacobi, J. Hadamard,
H. Poincar\'e, G. D. Birkhoff, M. Morse, L. Lyusternik and Schnirelmann
and others. Specially G. D. Birkhoff established the existence of at least
one closed geodesic on every Riemannian sphere $S^d$ with $d\ge 2$
(cf. \cite{Bir1}). Later L. Lyusternik and A. Fet proved the existence of
at least one closed geodesic on every compact Riemannian manifold (cf.
\cite{LyF1}). Such a variational proof works also for Finlser metrics
on compact manifolds and produces at least one closed geodesic on every
such manifold. An important breakthrough on this study is due to V. Bangert
\cite{Ban2} and J. Franks \cite{Fra1} around 1990, who proved the existence
of infinitely many closed geodesics on every Riemannian $2$-sphere (cf.
also \cite{Hin2} and \cite{Hin3} for new proofs of some parts of this result).

For irreversible Finsler manifolds, the closed geodesic problem is more
delicate as discovered by A. Katok via his famous example of 1973 which yields
some irreversible Finsler metrics on $S^d$ with precisely $2[(d+1)/2]$ distinct
prime closed geodesics (cf. \cite{Kat1} and \cite{Zil1}). In \cite{HWZ1} of
2003, H. Hofer, K. Wysocki and E. Zehnder proved that there exist either two
or infinitely many distinct prime closed geodesics on a Finsler $(S^2,F)$
provided that all the iterates of all closed geodesics are non-degenerate and
the stable and unstable manifolds of all hyperbolic closed geodesics intersect
transversally. In \cite{BaL1} of 2005 published in 2010, V. Bangert and Y. Long
proved that on every irreversible Finsler $S^2$ there exist always at least two
distinct prime closed geodesics (cf. also \cite{LoW2}).

Here recall that on a Finsler manifold $(M,F)$, a closed geodesic
$c:S^1=\R/\Z\to M$ is {\it prime}, if it is not a multiple covering (i.e., iteration)
of any other closed geodesic. Here the $m$-th iteration $c^m$ of $c$ is defined by
$c^m(t)=c(mt)$ for $m\in\N$. The inverse curve $c^{-1}$ of $c$ is defined by
$c^{-1}(t)=c(1-t)$ for $t\in S^1$. Two prime closed geodesics $c_1$ and $c_2$ on
a Finsler manifold $(M,F)$ (or Riemannian manifold $(M,g)$) are {\it distinct}
(or {\it geometrically distinct}), if they do not differ by an $S^1$-action
(or $O(2)$-action). We denote by $\CG(M,F)$ the set of all distinct closed geodesics
on $(M,F)$ for Finsler or Riemannian metric $F$ on $M$.

A long-standing conjecture on the closed geodesics is
\be  \;^{\#}\CG(M,g) = +\infty,  \lb{1.1}\ee
for every Riemannian metric $g$ on any compact manifold $M$ with $\dim M\ge 2$.
Correspondingly for Finsler manifolds, it is conjectured (cf. \cite{Lon6}) that
for each positive integer $n$ there exist positive integers $1\le p_n\le q_n$ with
$p_n\to +\infty$ as $n\to +\infty$ such that there holds
\be  \;^{\#}\CG(M,F) \in [p_n,\,q_n]\cup \{+\infty\},   \lb{1.2}\ee
for every Finsler metric $F$ on each compact manifold $M$ satisfying $\dim M=n$.

Note that by the results of \cite{Ban2} and \cite{Fra1} and the classification of
$2$-dimensional compact manifolds, the conjecture (\ref{1.1}) was proved when
$\dim M=2$. Similarly by the results of \cite{Kat1} and \cite{BaL1}, we have $p_2=2$.

In the study of the conjecture (\ref{1.1}), D. Gromoll and W. Meyer \cite{GrM1}
in 1969 proved the following result:

{\bf Theorem A.} (\cite{GrM1}) {\it On a compact Riemannian manifold
there exist infinitely many closed geodesics, if the free loop space of
this manifold has an unbounded sequence of Betti numbers.}

Stimulated by this result, M. Vigu\'e-Poirrier and D. Sullivan \cite{ViS1}
in 1976 proved:

{\bf Theorem B.} (\cite{ViS1}) {\it The free loop space of a
compact simply connected Riemannian manifold $M$ has no unbounded sequence
of Betti numbers if and only if the rational cohomology algebra of $M$
possess only one generator. }

Both of these two theorems were generalized to corresponding Finsler manifolds by
H. Matthias in 1980 (cf. \cite{Mat1}). Therefore based on these two theorems,
the most interesting manifolds in this multiplicity problem are those
compact simply connected manifolds satisfying
\be  H^*(M;\Q)\cong T_{d,h+1}(x)=\Q[x]/(x^{h+1}=0)   \lb{1.3}\ee
with a generator $x$ of degree $d\ge 2$ and hight $h+1\ge 2$. The most important
examples here are certainly spheres $S^d$ of dimension $d$.

Besides these results, when the dimension of a compact simply connected manifold
is greater than $2$, we are not aware of any multiplicity results on the existence
of at least two closed geodesics without pinching, generic, bumpy or other
conditions even on spheres (cf. \cite{Ano1}, \cite{Ban1}, \cite{Kli1}, \cite{BTZ1},
\cite{BTZ2}, \cite{DuL1}, \cite{DuL2}, \cite{LoW1}, \cite{Rad3}, \cite{Rad4},
\cite{Rad5}, \cite{Rad6}), except the Theorem C below proved recently in \cite{LoD1}
for the $3$-dimensional case and \cite{DuL3} for the $4$-dimensional case.

{\bf Theorem C.} (\cite{LoD1}, \cite{DuL3}) {\it There exist always at least two
distinct prime (geometrically distinct) closed geodesics for every irreversible
(or reversible, specially Riemannian) Finsler metric on every $3$ or
$4$-dimensional compact simply connected manifold. }

\smallskip

In this paper, we further generalize Theorem C to all compact simply
connected Finsler as well as Riemannian manifolds and prove the following
results.

{\bf Theorem 1.1.} {\it For every irreversible Finsler metric $F$ on
any compact simply connected manifold of dimension at least $2$, there exist
always at least two distinct prime closed geodesics.}

\smallskip

{\bf Theorem 1.2.} {\it For every reversible Finsler metric $F$ on
any compact simply connected manifold of dimension at least $2$, there exist
always at least two geometrically distinct closed geodesics. In particular,
it holds for every such Riemannian manifold.}

\smallskip

Next we briefly describe the main ideas in the proofs of Theorems 1.1 and 1.2.

In \cite{LoD1}, we classified all the closed geodesics into {\it rational} and
{\it irrational} two classes according to the basic normal form decomposition of
their linearized Poincar\'e maps as symplectic matrices introduced in \cite{Lon2}
in 1999. Then in \cite{LoD1}, we established periodicity of the Morse indices
and homological information of iterates of orientable rational closed geodesics
on any Finsler manifold $(M,F)$. Specially we proved
\be  i(c^{n+m}) = i(c^n) + i(c^m) + \ol{p}(c), \qquad \nu(c^{n+m})=\nu(c^m),
             \quad \forall\, m\in\N, \lb{1.4}\ee
where $n=n(c)$ is the analytical period of a prime closed geodesic $c$, cf. (\ref{4.1})
below, and $\ol{p}(c)$ is a constant depends only on the linearized Poincar\'e
map $P_c$ of $c$. We proved also a boundedness property of Morse indices in iterates of
every prime orientable rational closed geodesic $c$:
\be  i(c^m) + \nu(c^m) \le i(c^n) + \ol{p}(c) + \dim M -3,
             \quad \forall\,1\le m\le n-1.       \lb{1.5}\ee
If $(M,F)$ is a compact simply connected Finsler manifold and possesses only one
prime closed geodesic $c$, and if $c$ is rational, based on the properties (\ref{1.4})
and (\ref{1.5}) we established in \cite{LoD1} and \cite{DuL3} the following identity
\be B(d,h)(i(c^n) + \ol{p}(c)) + (-1)^{i(c^n)+\mu}\ka = \sum_{j=\mu-\ol{p}(c)+1}^{i(c^n)+\mu}(-1)^j b_j,
                 \lb{1.6}\ee
for some integer $\ka\ge 0$, where $\mu=\ol{p}(c)+\dim M -3$, and $B(d,h)$ depends
only on $d$ and $h$ and is given in Lemma 2.4 below, $b_j$s are Betti numbers of the
relative free loop spaces defined in Lemmas 2.5 and 2.6 below. Then using (\ref{1.6}),
and our computations on the precise sum of Betti numbers, we obtain a contradiction
and conclude that the only one prime closed geodesic $c$ on $M$ can not be rational.

Now in the current paper, our main idea is to generalize the above method on rational
closed geodesics to every closed geodesic on compact simply connected manifolds. Suppose
that there exists only one prime closed geodesic $c$ on a compact simply connected
Finsler manifold $(M,F)$. The new observations in the current paper are the following:

(i) When $c$ is irrational, suppose the basic normal form decomposition of the
linearized Poincar\'{e} map $P_c$ of $c$ contains $k$ irrational rotation matrices.
In this case, we can not hope the periodicity (\ref{1.4}) to be still true
anymore for the analytical period $n=n(c)$ and the constant $\ol{p}(c)$. But using the mod
one uniform distribution property of irrational numbers, we can still get a local version
of (\ref{1.4}), i.e., there exists a large enough even integer $T\in n(c)\N$ such that
for some integer $m_0=m_0(c)>1$ depending on $c$ only there holds
\be  i(c^{T+m}) = i(c^T) + i(c^m) + p(c), \quad \nu(c^{T+m})=\nu(c^m),
             \qquad \forall\,1\le m\le m_0, \lb{1.7}\ee
where $p(c)=\ol{p}(c)+2(A-k)$ for some integer $1\le A\le k$ depending on $P_c$ only.
We call such a property the {\it quasi-periodicity} of Morse indices of iterates
$c^m$.

(ii) Similarly for irrational $c$, we can not hope (\ref{1.5}) to be true for all
multiples of $n$. But using estimates on Morse indices of iterates of irrational closed
geodesics established in \cite{DuL3}, we can get also a similar version of (\ref{1.5}),
i.e., we can further choose the integer $T\in n(c)\N$ so that there holds
\be  i(c^m) + \nu(c^m) \le i(c^T) + p(c) + \dim M -3, \qquad
                 \forall\,1\le m\le T-1.   \lb{1.8}\ee

(iii) Now by computing out the alternating sum of the dimensions of all the critical
modules of $c^m$ with $1\le m \le T$, and then comparing with the Betti numbers of
the free loop space pairs on $M$ (i.e., $b_j$s below), we obtain the following version
of (\ref{1.6}) which holds at the iteration $T$: i.e., there exists an integer
$\ka\ge 0$ such that
\be B(d,h)(i(c^T) + p(c)) + (-1)^{i(c^T)+\mu}\ka = \sum_{j=\mu-p(c)+1}^{i(c^T)+\mu}(-1)^j b_j,
                 \lb{1.9}\ee
where $\mu=p(c)+\dim M -3$. Note that (\ref{1.7})-(\ref{1.9}) are automatically reduced
to (\ref{1.4})-(\ref{1.6}) when $c$ is rational.

(iv) Then the precise sum of Betti numbers on the right hand side of (\ref{1.9})
yields a contradiction, and shows that there must exist at least two distinct closed
geodesics.

Here we should point out that the identity (\ref{1.9}) (or (\ref{1.6})) is rather
different from the Morse inequalities, because the term $B(d,h)(i(c^T)+p(c))$ in
(\ref{1.9}) (or the corresponding term in (\ref{1.6})) represents the alternating
sum of dimensions of all local critical modules of $c^m$ with $1\le m\le T$
(or $1\le m\le n$), which is the alternating sum of all terms on or below the
$T$th (or $n$th) horizontal line in the Figure (\ref{5.54}) below, and is not the
alternating sum of Morse type numbers of $c^m$s with dimensions less than a fixed
integer, which is the alternating sum of all terms on the left of some fixed
vertical line in the Figure (\ref{5.54}) below. In fact in our case, firstly the
alternating sum of Morse type numbers with dimensions less than some fixed integer
in the Morse inequality may not be computable, because in general there may not
exist such a vertical line in the Figure (\ref{5.54}) below such that all
non-trivial critical modules of each iterate $c^m$ appears only on one side of
this vertical line. Secondly, even if it is computable, it is still not clear
whether the corresponding Morse inequalities may yield any contradiction.

Note that in his famous book \cite{Mor1}, M. Morse proved that for any given integer
$N>0$ the global homology of a $d$-dimensional ellipsoid ${\cal E}_d$ at all dimensions
less than $N$ can be produced by iterates of the $(d+1)$ main ellipses only, provided
${\cal E}_d$ is sufficiently close to the ball and all of its semi-axis are different.
His this example explains why the iterate $T$ in our proof should be sufficiently large
and carefully chosen.

\smallskip

For reader's conveniences, in Section 2 we briefly review some known results
on closed geodesics and Betti numbers of the $S^1$-invariant free loop space
of compact simply connected manifolds satisfying the condition (\ref{1.3}).
In Section 3 we briefly review basic normal form decompositions of symplectic
matrices and the precise index iteration formulae of symplectic paths established by
Y. Long in \cite{Lon2} and \cite{Lon3} together with the orientability of closed
geodesics. In Section 4, we establish the quasi-periodicity (\ref{1.7}) and the
boundedness estimate (\ref{1.8}) of iterated indices of closed geodesics. In Section 5,
using the index quasi-periodicity we prove some homological isomorphism theorems of
energy critical level pairs when there exists only one prime closed geodesic, and then
establish the identity (\ref{1.9}). In Section 6 we give proofs of Theorems 1.1 and 1.2.

In this paper, we denote by $\N$, $\N_0$, $\Z$, $\Q$, $\R$, and $\C$ the
sets of positive integers, non-negative integers, integers, rational numbers,
real numbers, and complex numbers respectively. We define the functions
$[a]=\max\{k\in\Z\,|\,k\le a\}$, $\{a\}=a-[a]$,
$E(a)=\min\{k\in\Z\,|\,k\ge a\}$ and $\vf(a)=E(a)-[a]$. Denote by $\,^{\#}A$ the
number of elements in a finite set $A$. When $S^1$ acts on a topological space $X$,
we denote by $\ol{X}$ the quotient space $X/S^1$. In this paper, we use only
singular homology modules with $\Q$-coefficients.

\setcounter{equation}{0}
\section{Critical point theory of closed geodesics}

\subsection{Critical modules for closed geodesics}

Let $M$ be a manifold with a Finsler metric $F$. Closed geodesics are critical
points of the energy functional
$E(\gamma)=\frac{1}{2}\int_{S^1}F(\gamma(t),\dot{\gamma}(t))^2dt$
on the Hilbert manifold $\Lambda M$ of $H^1$-maps from $S^1$ to $M$. An
$S^1$-action is defined by $(s\cdot\gamma)(t)=\gamma(t+s)$ for all
$\gamma\in\Lm M$ and $s, t\in S^1$. The index form of the functional
$E$ is well defined along any closed geodesic $c$ on $M$, which we
denote by $E''(c)$. As usual, denote by $i(c)$ and $\nu(c)$ the
Morse index and nullity of $E$ at $c$. For a closed geodesic $c$,
denote by $c^m$ the $m$-fold iteration of $c$ and
$\Lm(c^m)=\{\ga\in\Lm M\,|\, E(\ga)<E(c^m)\}$. Recall that
respectively the {\it mean index} $\hat{i}(c)$ and the $S^1$-{\it
critical modules} of $c^m$ are defined by
\be \hat{i}(c)=\lim_{m\rightarrow\infty}\frac{i(c^m)}{m},
   \quad \ol{C}_*(E,c^m) = H_*\left((\Lm(c^m)\cup S^1\cdot c^m)/S^1,\Lm(c^m)/S^1\right).
        \lb{2.1}\ee

If $c$ has multiplicity $m$, then the subgroup $\Z_m=\{\frac{n}{m}:0\le n<m\}$
of $S^1$ acts on $\ol{C}_k(E,c)$. As on page 59 of \cite{Rad2}, for $m\ge1$, let
$H_*(X,A)^{\pm \Z_m}=\{[\xi]\in H_*(X,A):T_*[\xi]=\pm\xi\}$, where $T$ is a
generator of the $\Z_m$ action. On $S^1$-critical modules of $c^m$, the following
lemma holds:

\smallskip

{\bf Lemma 2.1.} (cf. \cite{Rad2}, \cite{BaL1}, \cite{LoD1}) {\it Suppose $c$ is
a prime closed geodesic on a Finsler manifold $M$. Then
there exist two sets $U_{c^m}^-$ and $N_{c^m}$, the so-called local negative
disk and the local characteristic manifold at $c^m$ respectively,
such that $\nu(c^m)=\dim N_{c^m}$ and
\bea \overline{C}_q( E,c^m)
&\equiv& H_q\left((\Lm(c^m)\cup S^1\cdot c^m)/S^1, \Lm(c^m)/S^1\right)\nn\\
&=& \left(H_{i(c^m)}(U_{c^m}^-\cup\{c^m\},U_{c^m}^-)
    \otimes H_{q-i(c^m)}(N_{c^m}^-\cup\{c^m\},N_{c^m}^-)\right)^{+\Z_m}, \nn\eea

(i) When $\nu(c^m)=0$, there holds
$$ \overline{C}_q( E,c^m) = \left\{\matrix{
     \Q, &\quad {\it if}\;\; i(c^m)=i(c)\,(\mod 2)\;\;{\it and}\;\;
                   q=i(c^m),\;  \cr
     0, &\quad {\it otherwise}, \cr}\right.  $$

(ii) When $\nu(c^m)>0$, let $\ep(c^m)=(-1)^{i(c^m)-i(c)}$, then
there holds
$$ \overline{C}_q( E,c^m)=H_{q-i(c^m)}(N_{c^m}^-\cup\{c^m\},N_{c^m}^-)^{\ep(c^m)\Z_m}. $$}

Let
\be  k_j(c^m) \equiv \dim\, H_j( N_{c^m}^-\cup\{c^m\},N_{c^m}^-), \quad
     k_j^{\pm 1}(c^m) \equiv \dim\, H_j(N_{c^m}^-\cup\{c^m\},N_{c^m}^- )^{\pm\Z_m}.
        \lb{2.2}\ee
Then we have

\smallskip

{\bf Lemma 2.2.} (cf. \cite{Rad2}, \cite{BaL1}, \cite{LoD1}) {\it Let $c$
be a closed geodesic on a Finsler manifold $M$.

(i) There hold $0\le k_j^{\pm 1}(c^m) \le k_j(c^m)$ for $m\ge 1$ and
$j\in\Z$, $k_j(c^m)=0$ whenever $j\not\in [0,\nu(c^m)]$ and
$k_0(c^m)+k_{\nu(c^m)}(c^m)\le1$. If $k_0(c^m)+k_{\nu(c^m)}(c^m)=1$,
then $k_j(c^m)=0$ when $j\in(0,\nu(c^m))$.

(ii) For any $m\in\N$, there hold $k_0^{+1}(c^m) = k_0(c^m)$ and
$k_0^{-1}(c^m) = 0$. In particular, if $c^m$ is non-degenerate,
there hold $k_0^{+1}(c^m) = k_0(c^m)=1$, and $k_0^{-1}(c^m) =
k_j^{\pm 1}(c^m)=0$ for all $j\neq 0$.

(iii) Suppose for some integer $m=np\ge 2$ with $n$ and $p\in\N$ the
nullities satisfy $\nu(c^m)=\nu(c^n)$. Then there hold
$k_j(c^m)=k_j(c^n)$ and ${k}_j^{\pm 1}(c^m)={k}^{\pm 1}_j(c^n)$ for
any integer $j$.}

\medskip

Let $(M,F)$ be a compact and simply connected Finsler manifold with
finitely many prime closed geodesics. It is well known that for every
prime closed geodesic $c$ on $(M,F)$, there holds either
$\hat{i}(c)>0$ and then $i(c^m)\to +\infty$ as $m\to +\infty$, or
$\hat{i}(c)=0$ and then $i(c^m)=0$ for all $m\in\N$. Denote those
prime closed geodesics on $(M,F)$ with positive mean indices by
$\{c_j\}_{1\le j\le k}$. In \cite{Rad1} and \cite{Rad2}, Rademacher
established a celebrated mean index identity relating all the $c_j$s
with the global homology of $M$ (cf. Section 7, specially Satz 7.9
of \cite{Rad2}) for compact simply connected Finsler manifolds. A
refined version of this identity with precise coefficients was
proved in \cite{BaL1}, \cite{LoW1}, and \cite{LoD1}.

For each $m\in\N$, let $\ep=\ep(c^m)=(-1)^{i(c^m)-i(c)}$ and
\bea  K(c^m)
&\equiv&(k_0^\ep(c^m), k_1^\ep(c^m), \ldots, k_{2\dim M-2}^\ep(c^m))\nn\\
&=&(k_0^{\ep (c^m)}(c^m),  k_1^{\ep(c^m)}(c^m), \ldots,
     k_{\nu (c^m)}^{\ep (c^m)}(c^m), 0, \ldots, 0). \lb{2.3}\eea

{\bf Lemma 2.3.} (cf. Lemmas 7.1 and 7.2 of \cite{Rad2}, cf. also \cite{LoD1})
{\it Let $c$ be a prime closed geodesic on a compact Finsler manifold $(M,F)$.
Then there exists a minimal integer $N=N(c)\in\N$ such that $\nu(c^{m+N})=\nu(c^m)$,
$i(c^{m+N})-i(c^m)\in 2\Z$, and $K(c^{m+N})=K(c^m)$ for all $m\in\N$.}

\smallskip

{\bf Lemma 2.4.} (cf. \cite{Rad2}, \cite{BaL1}, \cite{LoW1}, \cite{LoD1})
{\it Let $(M,F)$ be a compact simply connected Finsler manifold with
$\,H^{\ast}(M,\Q)=T_{d,h+1}(x)$ for some integers $d\ge 2$ and $h\ge 1$.
Denote prime closed geodesics on $(M,F)$
with positive mean indices by $\{c_j\}_{1\le j\le k}$ for some $k\in\N$.
Then the following identity holds
\be \sum_{j=1}^k\frac{\hat{\chi}(c_j)}{\hat{i}(c_j)}=B(d,h)
=\left\{\matrix{
     -\frac{h(h+1)d}{2d(h+1)-4}, &\quad d\,{\rm even},\cr
     \frac{d+1}{2d-2}, &\quad d\,{\rm odd}, \cr}\right.   \lb{2.4}\ee
where $\dim M=hd$, $h=1$ when $M$ is a sphere $S^d$ of dimension $d$ and
\be \hat{\chi}(c) = \frac{1}{N(c)}
\sum_{0\le l_m\le \nu(c^m) \atop 1\le m\le N(c)}(-1)^{i(c^m)+l_m}k_{l_m}^{\ep(c^m)}(c^m)\;
                 \in \;\Q.        \lb{2.5}\ee}

\subsection{The structure of $H_*(\Lm M/S^1, \Lm^0 M/S^1;\Q)$}

Set $\ol{\Lm}^0=\ol{\Lambda}^0M =\{{\rm
constant\;point\;curves\;in\;}M\}\cong M$. Let $(X,Y)$ be a
space pair such that the Betti numbers $b_i=b_i(X,Y)=\dim
H_i(X,Y;\Q)$ are finite for all $i\in \Z$. As usual the {\it
Poincar\'e series} of $(X,Y)$ is defined by the formal power series
$P(X, Y)=\sum_{i=0}^{\infty}b_it^i$. We need the following version
of results on Betti numbers. The precise computations on
each Betti number in Lemma 2.6 and sums of Betti numbers in Lemmas 2.5
and 2.6 were given in \cite{LoD1} and \cite{DuL3}.

{\bf Lemma 2.5.} (cf. Theorem 2.4 and Remark 2.5 of \cite{Rad1}, Proposition 2.4
of \cite{LoD1}, Lemma 2.5 of \cite{DuL3}) {\it Let $(S^d,F)$ be a
$d$-dimensional Finsler sphere.}

(i) {\it When $d$ is odd, the Betti numbers are given by
\bea b_j
&=& \rank H_j(\Lm S^d/S^1,\Lm^0 S^d/S^1;\Q)  \nn\\
&=& \left\{\matrix{
    2,&\quad {\it if}\quad j\in \K\equiv \{k(d-1)\,|\,2\le k\in\N\},  \cr
    1,&\quad {\it if}\quad j\in \{d-1+2k\,|\,k\in\N_0\}\bs\K,  \cr
    0 &\quad {\it otherwise}. \cr}\right. \lb{2.6}\eea
For any $k\in \N$ and $k\ge d-1$, there holds }
\bea \sum_{j=0}^k(-1)^jb_j
&=& \sum_{0\le 2j\le k}b_{2j}   \nn\\
&=& \frac{k(d+1)}{2(d-1)} - \frac{d-1}{2} - \ep_{d,1}(k)  \nn\\
&\le& \frac{k(d+1)}{2(d-1)} - \frac{d-1}{2}.   \lb{2.7}\eea
where $\ep_{d,1}(k) = \{\frac{k}{d-1}\} + \{\frac{k}{2}\}\in [0,\frac{3}{2}-\frac{1}{2(d-1)})$.

(ii) {\it When $d$ is even, the Betti numbers are given by
\bea b_j
&=& \rank H_j(\Lm S^d/S^1,\Lm^0 S^d/S^1;\Q)  \nn\\
&=& \left\{\matrix{
    2,&\quad {\it if}\quad j\in \K\equiv \{k(d-1)\,|\,3\le k\in (2\N+1)\},  \cr
    1,&\quad {\it if}\quad j\in \{d-1+2k\,|\,k\in\N_0\}\bs\K,  \cr
    0 &\quad {\it otherwise}. \cr}\right.    \lb{2.8}\eea
For any $k\in \N$ and $k\ge d-1$, there holds
\bea -\sum_{j=0}^k(-1)^jb_j
= \sum_{0\le 2j-1\le k}b_{2j-1}
\le \frac{k d}{2(d-1)} - \frac{d-2}{2}.   \lb{2.9}\eea}

\medskip

For a compact and simply connected Finsler manifold $M$ with
$H^*(M;\Q)\cong T_{d,h+1}(x)$, when $d$ is odd, then $x^2=0$ and
$h=1$ in $T_{d,h+1}(x)$. Thus $M$ is rationally homotopy equivalent to
$S^d$ (cf. Remark 2.5 of \cite{Rad1}). Therefore, next we only consider
the case when $d$ is even.

\medskip

{\bf Lemma 2.6.} (cf. Theorem 2.4 of \cite{Rad1}, Lemma 2.6 of \cite{DuL3})
{\it Let $M$ be a compact simply connected manifold with $H^*(M;\Q)\cong T_{d,h+1}(x)$
for some integer $h\ge 1$ and even integer $d\ge 2$. Let $D=d(h+1)-2$ and
\bea \Om(d,h) = \{k\in 2\N-1&\,|\,& iD\le k-(d-1)=iD+jd\le iD+(h-1)d\;  \nn\\
         && \mbox{for some}\;i\in\N\;\mbox{and}\;j\in [1,h-1]\}. \lb{2.10}\eea
Then the Betti numbers of the free loop space of $M$ defined by
$b_q = \rank H_q(\Lm M/S^1,\Lm^0 M/S^1;\Q)$ for $q\in\Z$ are given by
\be b_q = \left\{\matrix{
    0, & \quad \mbox{if}\ q\ \mbox{is even or}\ q\le d-2,  \cr
    [\frac{q-(d-1)}{d}]+1, & \quad \mbox{if}\ q\in 2\N-1\;\mbox{and}\;d-1\le q < d-1+(h-1)d, \cr
    h+1, & \quad \mbox{if}\ q\in \Om(d,h), \cr
    h, & \quad \mbox{otherwise}. \cr}\right.
               \lb{2.11}\ee
For every integer $k\ge d-1+(h-1)d=hd-1$, we have
\bea \sum_{q=0}^kb_q
&=& \frac{h(h+1)d}{2D}(k-(d-1)) - \frac{h(h-1)d}{4} + 1 + \ep_{d,h}(k) \nn\\
&<&  h(\frac{D}{2}+1)\frac{k-(d-1)}{D} - \frac{h(h-1)d}{4} + 2,   \lb{2.12}\eea
where
\bea \ep_{d,h}(k)
&=& \{\frac{D}{hd}\{\frac{k-(d-1)}{D}\}\} - (\frac{2}{d}+\frac{d-2}{hd})\{\frac{k-(d-1)}{D}\}   \nn\\
&&\quad - h\{\frac{D}{2}\{\frac{k-(d-1)}{D}\}\} - \{\frac{D}{d}\{\frac{k-(d-1)}{D}\}\},  \lb{2.13}\eea
and there hold $\ep_{d,h}(k)\in (-(h+2),1)$ and $\ep_{d,1}(k)\in (-2,0]$ for all integer $k\ge d-1$. }

\setcounter{equation}{0}
\section{A review of the precise index iteration formulae for symplectic paths}

For $d\in\N$ and $\tau>0$, denote by $\Sp(2d)$ the symplectic group whose elements are $2d\times 2d$
real symplectic matrices and let
$$  \P_{\tau}(2d)=\{\ga\in C([0,\tau],\Sp(2d))\,|\,\ga(0)=I\}. $$
An index function theory $(i_{\om}(\ga),\nu_{\om}(\ga))$ for every symplectic path
$\ga\in\P_{\tau}(2d)$ parametrized by $\om\in \U=\{z\in\C\,|,|z|=1\}$ was introduced
by Y. Long in \cite{Lon2} of 1999. This index function theory is based on the
Maslov-type index theory $(i_1(\ga),\nu_1(\ga))$ for symplectic paths in $\P_{\tau}(2d)$
established by C. Conley, E. Zehnder in \cite{CoZ1} of 1984, Y. Long and E. Zehnder in
\cite{LZe1} of 1990, and Y. Long in \cite{Lon1} of 1990 (cf. \cite{Lon5}). In \cite{Lon2},
Y. Long established also the basic normal form decomposition of symplectic matrices.
Based on this result he further established the precise iteration formulae of indices
of symplectic paths in \cite{Lon3} of 2000. These results form the basis of our study
on the Morse indices and homological properties of iterates of closed geodesics. Here
we briefly review these results.

As in \cite{Lon5}, denote by
\bea
N_1(\lm, a) &=& \left(\matrix{\lm & a\cr
                                0 & \lm\cr}\right), \qquad {\rm for\;}\lm=\pm 1, \; a\in\R, \lb{3.1}\\
H(b) &=& \left(\matrix{b & 0\cr
                      0 & b^{-1}\cr}\right), \qquad {\rm for\;}b\in\R\bs\{0, \pm 1\}, \lb{3.2}\\
R(\th) &=& \left(\matrix{\cos\th & -\sin\th \cr
                           \sin\th & \cos\th\cr}\right), \qquad {\rm for\;}\th\in (0,\pi)\cup (\pi,2\pi), \lb{3.3}\\
N_2(e^{\th\sqrt{-1}}, B) &=& \left(\matrix{ R(\th) & B \cr
                  0 & R(\th)\cr}\right), \qquad {\rm for\;}\th\in (0,\pi)\cup (\pi,2\pi)\;\; {\rm and}\; \nn\\
        &&B=\left(\matrix{b_1 & b_2\cr
                                  b_3 & b_4\cr}\right)\quad {\rm with}\; b_j\in\R, \;\;
                                         {\rm and}\;\; b_2\not= b_3. \lb{3.4}\eea
Here $N_2(e^{\th\sqrt{-1}}, B)$ is non-trivial if $(b_2-b_3)\sin\theta<0$, and trivial
if $(b_2-b_3)\sin\theta>0$. In \cite{Lon2}-\cite{Lon4}, these matrices are called
{\it basic normal forms} of symplectic matrices.

As in \cite{Lon5}, given any two real matrices of the square block form
$$ M_1=\left(\matrix{A_1 & B_1\cr C_1 & D_1\cr}\right)_{2i\times 2i},\qquad
   M_2=\left(\matrix{A_2 & B_2\cr C_2 & D_2\cr}\right)_{2j\times 2j},$$
the $\diamond$-sum (direct sum) of $M_1$ and $M_2$ is defined by the
$2(i+j)\times2(i+j)$ matrix
$$ M_1\diamond M_2=\left(\matrix{A_1 & 0 & B_1 & 0 \cr
                                   0 & A_2 & 0& B_2\cr
                                   C_1 & 0 & D_1 & 0 \cr
                                   0 & C_2 & 0 & D_2}\right). $$

{\bf Definition 3.1.} (cf. \cite{Lon3} and \cite{Lon5}) {\it For every
$P\in\Sp(2d)$, the homotopy set $\Omega(P)$ of $P$ in $\Sp(2d)$ is defined by
$$ \Om(P)=\{N\in\Sp(2d)\,|\,\sg(N)\cap\U=\sg(P)\cap\U\equiv\Gamma\;\mbox{and}
                    \;\nu_{\om}(N)=\nu_{\om}(P)\, \forall\om\in\Gamma\}, $$
where $\sg(P)$ denotes the spectrum of $P$,
$\nu_{\om}(P)\equiv\dim_{\C}\ker_{\C}(P-\om I)$ for all $\om\in\U$.
The homotopy component $\Om^0(P)$ of $P$ in $\Sp(2d)$ is defined by
the path connected component of $\Om(P)$ containing $P$ (cf. page 38 of
\cite{Lon5}). }

Note that $\Om^0(P)$ defines an equivalent relation among symplectic
matrices. Specially we call two matrices $N$ and $P\in\Sp(2d)$ {\it homotopic},
if $N\in\Om^0(P)$, and in this case we write $N\approx P$.

Then the following decomposition theorem is proved in \cite{Lon2}
and \cite{Lon3}

\medskip

{\bf Theorem 3.2.} (cf. Theorem 7.8 of \cite{Lon2}, Lemma 2.3.5 and
Theorem 1.8.10 of \cite{Lon5}) {\it For every $P\in\Sp(2d)$, there
exists a continuous path $f\in\Om^0(P)$ such that $f(0)=P$ and
\bea f(1)
&=& N_1(1,1)^{\dm p_-}\,\dm\,I_{2p_0}\,\dm\,N_1(1,-1)^{\dm p_+}  \nn\\
 &&\dm\,N_1(-1,1)^{\dm q_-}\,\dm\,(-I_{2q_0})\,\dm\,N_1(-1,-1)^{\dm q_+} \nn\\
 &&\dm\,R(\th_1)\,\dm\,\cdots\,\dm\,R(\th_k)\,\dm\,R(\th_{k+1})\,\dm\,\cdots\,\dm\,R(\th_r)  \nn\\
 &&\dm\,N_2(e^{\aa_{1}\sqrt{-1}},A_{1})\,\dm\,\cdots\,\dm\,N_2(e^{\aa_{k_{\ast}}\sqrt{-1}},A_{k_{\ast}})  \nn\\
 &&\qquad\qquad \,\dm\,N_2(e^{\aa_{k_{\ast}+1}\sqrt{-1}},A_{k_{\ast}+1})\,\dm\,\cdots
  \,\dm\,N_2(e^{\aa_{r_{\ast}}\sqrt{-1}},A_{r_{\ast}}) \nn\\
 &&\dm\,N_2(e^{\bb_{1}\sqrt{-1}},B_{1})\,\dm\,\cdots\,\dm\,N_2(e^{\bb_{k_{0}}\sqrt{-1}},B_{k_{0}})  \nn\\
 &&\qquad\qquad \,\dm\,N_2(e^{\bb_{k_0+1}\sqrt{-1}},B_{k_0+1})\,\dm\,\cdots\,
              \dm\,N_2(e^{\bb_{r_{0}}\sqrt{-1}},B_{r_{0}}) \nn\\
 &&\dm\,H(2)^{\dm h_+}\,\dm\,H(-2)^{\dm h_-}, \lb{3.5}\eea
where $\frac{\th_{j}}{2\pi}\not\in\Q$ for $1\le j\le k$ and
$\frac{\th_{j}}{2\pi}\in\Q$ for $k+1\le j\le r$;
$N_2(e^{\aa_{j}\sqrt{-1}},A_{j})$'s are nontrivial basic normal
forms with $\frac{\aa_{j}}{2\pi}\not\in\Q$ for $1\le j\le k_{\ast}$
and $\frac{\aa_{j}}{2\pi}\in\Q$ for $k_{\ast}+1\le j\le r_{\ast}$;
$N_2(e^{\bb_{j}\sqrt{-1}},B_{j})$'s are trivial basic normal forms
with $\frac{\bb_{j}}{2\pi}\not\in\Q$ for $1\le j\le k_0$ and
$\frac{\bb_{j}}{2\pi}\in\Q$ for $k_0+1\le j\le r_0$; $p_-=p_-(P)$,
$p_0=p_0(P)$, $p_+=p_+(P)$, $q_-=q_-(P)$, $q_0=q_0(P)$,
$q_+=q_+(P)$, $r=r(P)$, $k=k(P)$, $r_{j}=r_{j}(P)$, $k_j=k_j(P)$ with $j=\ast, 0$ and
$h_+=h_+(P)$ are nonnegative integers, and $h_-=h_-(P)\in \{0,1\}$;
$\th_j$, $\aa_j$, $\bb_j \in (0,\pi)\cup (\pi,2\pi)$; these integers
and real numbers are uniquely determined by $P$ and satisfy}
\be p_- + p_0 + p_+ + q_- + q_0 + q_+ + r + 2r_{\ast} + 2r_0 + h_- + h_+ = d. \lb{3.6}\ee

\medskip

Based on Theorem 3.2, the homotopy invariance and symplectic
additivity of the indices, the following precise iteration formula
was proved in \cite{Lon3}:

\medskip

{\bf Theorem 3.3.} (cf. \cite{Lon3}, Theorem 8.3.1 and Corollary
8.3.2 of \cite{Lon5}) {\it Let $\ga\in\P_{\tau}(2d)$. Denote the
basic normal form decomposition of $P\equiv \ga(\tau)$ by
(\ref{3.5}). Then we have
\bea i_1(\ga^m)
&=& m(i_1(\ga)+p_-+p_0-r ) + 2\sum_{j=1}^rE\left(\frac{m\th_j}{2\pi}\right) - r   \nn\\
&&  - p_- - p_0 - {{1+(-1)^m}\over 2}(q_0+q_+) \nn\\
&& + 2\sum_{j=k_{\ast}+1}^{r_{\ast}}\vf\left(\frac{m\aa_j}{2\pi}\right) - 2(r_{\ast}-k_{\ast}), \lb{3.7}\\
\nu_1(\ga^m) &=& \nu_1(\ga) + {{1+(-1)^m}\over 2}(q_-+2q_0+q_+) + 2\vs(m,\ga(\tau)),    \lb{3.8}\\
\hat{i}(\ga) &=& i_1(\ga) + p_- + p_0 - r +
\sum_{j=1}^r\frac{\th_j}{\pi},   \lb{3.9}\eea
where we denote by}
\bea
\vs(m,\ga(\tau)) &=& (r-k) - \sum_{j=k+1}^r\vf\left(\frac{m\th_j}{2\pi}\right)  \nn\\
&& + (r_{\ast}-k_{\ast}) -
\sum_{j=k_{\ast}+1}^{r_{\ast}}\vf\left(\frac{m\aa_j}{2\pi}\right)
             + (r_0-k_0) - \sum_{j=k_0+1}^{r_0}\vf\left(\frac{m\bb_j}{2\pi}\right).
             \lb{3.10}
\eea

Let
\be \M\equiv\{N_1(1, b_1),\,b_1=0,1; \;\;N_1(-1,b_2),\,b_2=0,\pm1;\;\;R(\th),
            \,\th\in(0,\pi)\cup(\pi,2\pi);\,H(-2)\}.               \lb{3.11}\ee
By Theorems 8.1.4-8.1.7 and 8.2.1-8.2.4 on pp179-187 of \cite{Lon5}, we have specially

\medskip

{\bf Proposition 3.4.} {\it Every path $\ga\in\P_{\tau}(2)$ with end
matrix homotopic to some matrix in $\M$ must have odd index
$i_1(\ga)$. Paths $\xi\in\P_{\tau}(2)$ ending at $N_1(1,-1)$ or $H(2)$ and
$\eta\in\P_{\tau}(4)$ with end matrix homotopic to $N_2(\om,B)$ must have
even indices $i_1(\xi)$ and $i_1(\eta)$.}

\medskip

The relation between the Morse indices of closed geodesics on Finsler manifolds
and the above Maslov-type index theory for symplectic paths was studied
by C. Liu and Y. Long in \cite{LLo1} and C. Liu in \cite{Liu1}. Specially
we have

\medskip

{\bf Proposition 3.5.} (Theorem 1.1 and Remark 4.2 of \cite{Liu1}, cf. also
Theorem 1.1 of \cite{LLo1}, ) {\it For any closed geodesic $c$ on a Finsler
manifold $(M,F)$ with $d=\dim M<+\infty$, denote its linearized Poincar\'e
map by $P_c$. Then there exists a path $\ga\in C([0,1],\Sp(2d-2))$
satisfying $\ga(0)=I$, $\ga(1)=P_c$, and}
\bea
(i(c),\nu(c)) &=& (i_1(\ga), \nu_1(\ga)), \qquad {\it if\;}c\;\;{\it is\;orientable},
                       \lb{3.12}\\
(i(c),\nu(c)) &=& (i_{-1}(\ga), \nu_{-1}(\ga)), \qquad {\it if\;}c\;\;
      {\it is\;unorientable\;and\;}d\;{\it is\;even}. \lb{3.13}\eea

By this result, the above index iteration formulae (Theorem 3.3) can be applied to
every orientable closed geodesic on Finsler and Riemannian manifolds. For unorientable
closed geodesics, one can get a similar iteration formulae using results in \cite{Lon5}.

\medskip

{\bf Remark 3.6.} Note that every closed geodesic $c$ on a simply connected Finsler
manifold is always orientable and thus Theorem 3.3 can be applied to get $i(c^m)$
directly (cf. Section 2.1-Appendix on pages 136-141 of \cite{Kli1}). In this
paper we are interested in orientable closed geodesics.

Next we need the following results from \cite{DuL3}.

\smallskip

{\bf Proposition 3.7.} (Corollary 3.19 of \cite{DuL3}) {\it Let
$v=(v_1,\ldots, v_k)\in (\R\bs\Q)^k$. Then there exists
an integer $A$ satisfying $[(k+1)/2]\le A\le k$ and a subset $P$ of
$\{1, \ldots, k\}$ containing $A$ integers, such that for any integer
$n\in\N$ and any small $\ep>0$ there exist infinitely many even integers
$T_1$ and $T_2\in n\N$ satisfying respectively }
\bea
&& \left\{\matrix{ \{T_1v_i\} > 1-\ep, &\qquad& {\it for}\;\;i\in P, \cr
                     \{T_1v_j\} < \ep, &\qquad& {\it for}\;\; j\in \{1,\ldots,k\}\bs P, \cr}\right.  \lb{3.14}\\
&{\it or}\quad & \left\{\matrix{ \{T_2v_i\} < \ep, &\qquad& {\it for}\;\;i\in P, \cr
                     \{T_2v_j\} > 1- \ep, &\qquad& {\it for}\;\; j\in \{1,\ldots,k\}\bs P. \cr}\right.
                      \lb{3.15}\eea

\smallskip

{\bf Theorem 3.8.} (Theorem 3.21 of \cite{DuL3}) ({\bf Quasi-monotonicity of index growth
for closed geodesics}) {\it Let $c$ be an orientable closed geodesic with mean index $\hat{i}(c)>0$ on
a Finsler manifold $(M,F)$ of dimension $d\ge 2$. Denote the basic
normal form decomposition of the linearized Poincar\'e map $P_c$ of $c$ by (\ref{3.5}). Then
there exist an integer $A$ with $[(k+1)/2]\le A\le k$ and a subset $P$ of integers
$\{1, \ldots, k\}$ with $A$ integers such that for any $\ep\in (0,1/4)$ there exist infinitely many
sufficiently large even integer $T\in n\N$ satisfying
\bea
\left\{\frac{T\th_j}{2\pi}\right\} &>& 1-\ep, \qquad {\it for}\;\; j\in P,   \lb{3.16}\\
\left\{\frac{T\th_j}{2\pi}\right\} &<& \ep, \qquad {\it for}\;\; j\in \{1,\ldots,k\}\bs P.  \lb{3.17}\eea
Consequently we have
\bea
i(c^m)-i(c^T) &\ge& K_1 \equiv \lm + (q_0+q_+) + 2(r-k) + 2(r_{\ast}-k_{\ast}) + 2A,
          \quad \forall\, m\ge T+1, \lb{3.18}\\
i(c^T)-i(c^m) &\ge& K_2 \equiv \lm - (q_0+q_+) + 2k - 2(r_{\ast}-k_{\ast}) - 2A,
          \quad \forall\, 1\le m\le T-1, \lb{3.19}\eea
where $\lm = i(c)+p_-+p_0-r$, the integers $p_-$, $p_0$, $q_0$, $q_+$, $r$, $k$, $r_{\ast}$
and $k_{\ast}$ are defined in (\ref{3.5}). }

\setcounter{equation}{0}
\section{Properties of Morse indices of iterates of closed geodesics}

Let $(M,F)$ be a Finsler manifold of dimension $d$. As in \cite{LoD1}, a matrix
$P\in\Sp(2d-2)$ is {\it rational}, if no basic normal form in (\ref{3.5}) of $P$
is of the form $R(\th)$ with $\th/\pi\in\R\bs\Q$, and is {\it irrational}, otherwise.
Let $c$ be a closed geodesic on $(M,F)$ whose linearized Poincar\'e map is
denoted by $P_c$ and then $P_c\in\Sp(2d-2)$. The closed geodesic $c$
is {\it rational}, {\it irrational}, if so is $P_c$. The {\it analytical period}
$n(c)$ of $c$ is defined by
\be n(c) = \min\{j\in\N\,|\,\nu(c^j)=\max_{m\ge 1}\nu(c^m),\;\;
                  i(c^{m+j})-i(c^{m})\in 2\Z, \;\;\forall\,m\in\N\}. \lb{4.1}\ee

One of the most important properties of $n=n(c)$ is
\be   n(c) = N(c),  \lb{4.2}\ee
where $N(c)$ is defined by Lemma 2.3. This was proved by Lemma 3.10 of \cite{DuL3}.

Next we need

{\bf Definition 4.1.} {\it For any closed geodesic $c$ with $\hat{i}(c)>0$ on
a Finsler manifold $(M,F)$ of dimension $d$, Denote the
basic normal form decomposition of the linearized Poincar\'e
map $P_c$ of $c$ by (\ref{3.5}). We define $m_0=m_0(c)$ as follows
$$ m_0(c)=\min\{m\in\N\,|\,i(c^{j+m})\ge d+4k,\,\forall\,j\ge1\}. $$
Here $k$ is defined in (\ref{3.5}). Note that $\hat{i}(c)>0$ implies that
$i(c^m)\rightarrow +\infty$ as $m\rightarrow +\infty$. Thus the integer $m_0$
is well-defined.}

\medskip

{\bf Lemma 4.2.} {\it Let $c$ be an orientable prime closed geodesic on a Finsler
manifold $M=(M,F)$ with $\dim M<+\infty$. Then for every $m\in\N$, we have
\be i(c^{2m}) = i(c^2) \;\; \mod \; 2, \quad i(c^{2m+1}) = i(c) \;\;\mod\;2. \lb{4.3}\ee
For any two positive integers $q|p$, we have}
\be  i(c^p) \ge i(c^q) \quad {\it and} \quad  \nu(c^p) \ge \nu(c^q).  \lb{4.4}\ee

{\bf Proof.} (\ref{4.3}) follows from (\ref{3.7}) of Theorem 3.3 immediately.
The two inequalities in (\ref{4.4}) follow from the Bott formulae (cf. Theorem 1 and
its corollary on pages 177-178 of \cite{Bot1}) immediately. In fact using notations
of \cite{Lon5} we have
$$ i(c^p) = \sum_{\om^p=1}i_{\om}(c) = \sum_{\om^q=1}i_{\om}(c) + \sum_{\om^p=1,\, \om^q\not=1}i_{\om}(c)
     \ge \sum_{\om^q=1}i_{\om}(c) = i(c^q), $$
and
$$  \nu(c^p) = \sum_{\om^p=1}\nu_{\om}(c)
= \sum_{\om^q=1}\nu_{\om}(c) + \sum_{\om^p=1,\, \om^q\not=1}\nu_{\om}(c)
\ge \sum_{\om^q=1}\nu_{\om}(c) = \nu(c^q).  $$
Here that both $i_{\om}(c)$ and $\nu_{\om}(c)$ are non-negative integers for
any $\om\in\U$ follow from Proposition 1.3 of \cite{Bot1}. The proof is complete. \hfill\hb

\medskip

The following theorem gives some precise index properties of iterates of closed
geodesics, which is a crucial step in the proofs of Theorems 1.1 and 1.2, and
generalizes Theorem 3.7 in \cite{LoD1} for rational closed geodesics to the
non-rational closed geodesics.

\medskip

{\bf Theorem 4.3.} ({\bf Index quasi-periodicity of closed geodesics})
{\it Let $c$ be an orientable closed geodesic with $\hat{i}(c)>0$ on a
Finsler manifold of dimension $d$. Denote the basic normal form decomposition
of the linearized Poincar\'e map $P_c$ of $c$ by (\ref{3.5}). Let $n=n(c)$ be
the analytical period of $c$.

Then when $k\ge 1$, there exist an integer $A$ with $[(k+1)/2]\le A\le k$ and a subset
$P$ of integers $\{1, \ldots, k\}$ with $A$ integers such that for any given
integer $m_0\in\N$ and any small $\ep>0$ there exists a sufficiently large even
integer $T\in n\N$ satisfying
\bea
\left\{\frac{T\th_j}{2\pi}\right\} &>& 1-\ep, \qquad {\it for}\;\; j\in P,   \lb{4.5}\\
\left\{\frac{T\th_j}{2\pi}\right\} &<& \ep, \qquad {\it for}\;\; j\in \{1,\ldots,k\}\bs P.  \lb{4.6}\eea
Note that both of (\ref{4.5}) and (\ref{4.6}) should be omitted when $k=0$, and the following
conclusions (A)-(D) still hold by Theorem 3.7 of \cite{LoD1}.

Let
\be p(c) \equiv p_- + p_0 + q_0 + q_+ +2r_{\ast}-2k_{\ast}+r+2A-2k \ge 0. \lb{4.7}\ee
Then the following conclusions hold always:

(A) (Quasi-periodicity) For any $1\le m\le m_0$, there hold
\bea
i(c^{m+T}) &=& i(c^m) + i(c^{T}) + p(c), \lb{4.8}\\
\nu(c^{m+T}) &=& \nu(c^m).  \lb{4.9}\eea

(B) (Relative parity) There holds
\be  i(c^T) = p(c) \quad (\mod \, 2).   \lb{4.10}\ee

(C) (Nullity-periodicity) There holds
\be  \nu(c^n)=\nu(c^T) \le p(c) + d-1-2A.    \lb{4.11}\ee

(D) (Period-mean index) If $\hat{i}(c)>0$ is a rational number, there holds
\be  T\hat{i}(c) = i(c^T) + p(c).  \lb{4.12}\ee
If $\hat{i}(c)>0$ is irrational, then for any small $\tau>0$ we can further
require the above chosen $T\in n\N$ to be even larger to satisfy }
\be  |T\hat{i}(c) - (i(c^T) + p(c))| < \tau.  \lb{4.13}\ee

\medskip

{\bf Proof.} Note that by the definitions of $E(\cdot)$ and $\vf(\cdot)$ there hold
\bea
E(z+b)=z+E(b), &&  \vf(z+b)=\vf(b)\;\;\;{\rm for}\;\;k\in\Z,\;\; b\not\in\Z, \lb{4.14}\\
E(b)+E(-b) =1, && {\rm for}\;\;b\in (0,1), \lb{4.15}\\
E(a)+E(-a) - \vf(a) = 0, && \vf(a)=\vf(-a) \quad \forall \; a\in\R. \lb{4.16}\eea

Let $n=n(c)$ be the analytical period of $c$. For the integer $m_0$ given in the
assumption of the theorem, we specially set
\be \ep = \min\left\{\{\frac{m\th_j}{2\pi}\},\,1-\{\frac{m\th_j}{2\pi}\}\,|\,1\le m\le m_0;\,1\le j\le k\right\}.
          \lb{4.17}\ee

Note that, when $k\ge 1$, we fix an even integer $T\in n\N$ obtained from Proposition 3.7
satisfying (\ref{4.5}) and (\ref{4.6}) for this $\ep>0$.

We use short hand notations as in (\ref{3.5}) and carry out the proof in several steps.

\medskip

{\bf Step 1.} {\it Proof of the quasi-periodicity (A). }

\medskip

By (\ref{3.7}) of Theorem 3.3 for any $m\in\N$ we obtain
\bea i(c^{m+T})
&=& (m+T)(i(c)+p_-+p_0-r ) + 2\sum_{j=1}^rE(\frac{(m+T)\th_j}{2\pi}) - r   \nn\\
&&  - p_- - p_0 - {{1+(-1)^m}\over 2}(q_0+q_+) + 2\sum_{j=k_{\ast}+1}^{r_{\ast}}\vf(\frac{m\aa_j}{2\pi})
      - 2(r_{\ast}-k_{\ast}) \nn\\
&=&  i(c^m)+i(c^T)  + (r + p_- + p_0 + q_0 + q_+ +2r_{\ast}-2k_{\ast})   \nn\\
&&  \qquad+ 2\sum_{j=1}^rE(\frac{(m+T)\th_j}{2\pi})-2\sum_{j=1}^rE(\frac{m\th_j}{2\pi})
      - 2\sum_{j=1}^rE(\frac{T\th_j}{2\pi}). \lb{4.18}\eea
where we have used the evenness of $T$ and the fact $T\in n\N$. Note that
\bea \frac{1}{2}\Theta(m,T)
&\equiv& \sum_{j=1}^rE(\frac{(m+T)\th_j}{2\pi})-\sum_{j=1}^rE(\frac{m\th_j}{2\pi})
                 -\sum_{j=1}^rE(\frac{T\th_j}{2\pi})  \nn\\
&=& \sum_{j=1}^kE(\frac{(m+T)\th_j}{2\pi})-\sum_{j=1}^kE(\frac{m\th_j}{2\pi})
            -\sum_{j=}^kE(\frac{T\th_j}{2\pi}) \nn\\
&=& \sum_{j=1}^kE\left(\{\frac{m\th_j}{2\pi}\}+\{\frac{T\th_j}{2\pi}\}\right)
        -\sum_{j=1}^kE\left(\{\frac{m\th_j}{2\pi}\}\right)
        -\sum_{j=1}^kE\left(\{\frac{T\th_j}{2\pi}\}\right)  \nn\\
&=& \sum_{j=1}^kE\left(\{\frac{m\th_j}{2\pi}\}+\{\frac{T\th_j}{2\pi}\}\right)-2k  \nn\\
&=& \sum_{j=1}^AE\left(\{\frac{m\th_j}{2\pi}\}+\{\frac{T\th_j}{2\pi}\}\right)
    +\sum_{j=A+1}^kE\left(\{\frac{m\th_j}{2\pi}\}+\{\frac{T\th_j}{2\pi}\}\right)-2k.  \lb{4.19}\eea
So it follows from (\ref{4.5})-(\ref{4.6}) and (\ref{4.19}) that
\be \Theta(m,T)=2(A-k),\qquad \forall\, 1\le m\le m_0.  \lb{4.20}\ee
Together with (\ref{4.18}), it yields
\be  i(c^{m+T})= i(c^m)+i(c^T)  + r + p_- + p_0 + q_0 + q_+ + 2(r_{\ast}-k_{\ast})+2(A-k),\quad
       \forall\,1\le m\le m_0.   \lb{4.21}\ee

Thus (\ref{4.8}) holds. And (\ref{4.9}) follows from the definition of $T\in n\N$.

\medskip

{\bf Step 2.} {\it Proof of the relative parity (B).}

\medskip

By Theorem 3.3 and the definition (\ref{4.7}) of $p(c)$ we have
\bea i(c^T)-p(c)
&=& T(i(c)+p_-+p_0-r )+2\sum_{j=1}^r E(\frac{T\th_j}{2\pi})  \nn\\
& & - r - p_- - p_0 - {{1+(-1)^T}\over 2}(q_0+q_+)-2(r_{\ast}- k_{\ast})  \nn\\
& & - (p_- + p_0 + q_+ + q_0 +2r_{\ast}-2k_{\ast}+ r +2A-2k )  \nn\\
&=& T(i(c)+p_-+p_0-r ) + 2\sum_{j=1}^r E(\frac{T\th_j}{2\pi}) - 2r - 2p_- - 2p_0   \nn\\
& &  - {{3+(-1)^T}\over 2}(q_0+q_+) - 4(r_{\ast}-k_{\ast})-2(A-k).  \nn\eea
Because $T$ is even, it yields the relative parity (B).

\medskip

{\bf Step 3.} {\it Proof of the nullity-periodicity (C).}

\medskip

Because $\nu(c)=p_-+2p_0+p_+$, by Theorem 3.3 we have
\bea \nu(c^T)-p(c)&=&\nu(c^n)-p(c)\nn\\
&=& p_- + 2p_0 + p_+ + (q_-+2q_0+q_+) + 2(r-k+r_{\ast}-k_{\ast}+r_0-k_0) \nn\\
&& \quad - (p_0+p_-+q_0+q_++2r_{\ast}-2k_{\ast}+r+2A-2k)  \nn\\
&=& p_0 + p_+ + q_0 + q_- + r + 2r_0-2k_0-2A   \nn\\
&\le& d-1-2A.   \nn\eea
This yields (C).

\medskip

{\bf Step 4.} {\it Proof of the period-mean index (D). }

When $k=0$, (D) is proved in Theorem 3.7 of \cite{LoD1}. Now we consider the case
$k\ge 1$.

When $\hat{i}(c)=i(c)+p_-+p_0-r+\sum_{j=1}^r\th_j/\pi$ is a rational
number, we must have $k\ge 2$ and then $A\ge 1$ by Proposition 3.7.
Let $\sum_{j=1}^k\th_j/2\pi=q/p$ for some integers $0<p,q\in\N$ with $(p,q)=1$.
Further choose $0<\ep<\frac{1}{k}$ and an even $T\in np\N$ satisfying (\ref{4.5})
and (\ref{4.6}). Note that $\sum_{j=1}^k\{\frac{T\th_j}{2\pi}\}$ is an integer
because $T$ is an integer multiple of $p$.

If $A=k\ge 2$, by (\ref{4.5}) it yields a contradiction
\be \sum_{j=1}^k\{\frac{T\th_j}{2\pi}\}=\sum_{j=1}^A\{\frac{T\th_j}{2\pi}\}\in (A(1-\ep),A)\cap\Z=\emptyset.
               \lb{4.22}\ee

If $[\frac{k+1}{2}]\le A\le k-1$, by (\ref{4.5}) and (\ref{4.6}) we obtain
\be \sum_{j=1}^k\{\frac{T\th_j}{2\pi}\}\in (A(1-\ep),A+(k-A)\ep)\cap\Z=\{A\}.  \lb{4.23}\ee
Together with (\ref{3.7}) of Theorem 3.3 and the definition (\ref{4.7}) of $p(c)$, it
yields
\bea  T\hat{i}(c)
&=& T\left(i(c) + p_- + p_0 - r + \sum_{j=1}^r\frac{\th_j}{\pi}\right)  \nn\\
&=& i(c^T) + (r+p_-+p_0+2r_{\ast}-2k_{\ast})\nn\\
&& +\frac{1+(-1)^T}{2}(q_0+q_+)
  + 2\left(\sum_{j=1}^k\frac{T\th_j}{2\pi} -\sum_{j=1}^k E(\frac{T\th_j}{2\pi})\right)\nn\\
&=& i(c^T) + (r+p_-+p_0+2r_{\ast}-2k_{\ast}+q_0+q_+)   \nn\\
&& +2\left(\sum_{j=1}^k\{\frac{T\th_j}{2\pi}\} -\sum_{j=1}^k E(\{\frac{T\th_j}{2\pi}\})\right)\nn\\
&=& i(c^T)+(r+p_-+p_0+2r_{\ast}-2k_{\ast}+q_0+q_++2A-2k)\nn\\
&=& i(c^T) + p(c), \lb{4.24}\eea
where we have used the fact that $T\in np\N$ is even, and the
rationality of $\hat{i}(c)$ is used only to get the second last equality.

When $\hat{i}(c)>0$ is irrational, we further require that $\ep>0$ satisfies $2k\ep<\tau$.
Thus in this case (\ref{4.22}) and (\ref{4.23}) become
$$ \sum_{j=1}^k\{\frac{T\th_j}{2\pi}\}\in (A(1-\ep),A+k\ep)\subset (A-\tau/2, A+\tau/2). $$
Then from the third equality of (\ref{4.24}) we obtain
$$  |T\hat{i}(c) - (i(c^T) + p(c))| \le 2|\sum_{j=1}^k\{\frac{T\th_j}{2\pi}\}-A| < \tau. $$
i.e., (D) holds.

This completes the proof of Theorem 4.3. \hfill\hb

\medskip

{\bf Lemma 4.4.} {\it For any orientable closed geodesic $c$ with $\hat{i}(c)>0$ on a
Finsler manifold $(M,F)$ of dimension $d$, denote the basic normal form decomposition
of the linearized Poincar\'e map $P_c$ of $c$ by (\ref{3.5}). Then for any even $T\in n\N$
and $m_0=m_0(c)$ given by Definition 4.1, there holds
\bea i(c^{T+m_0+m})-i(c^T)\ge p(c)+d, \qquad\forall\, m\ge 1. \lb{4.25}\eea}

{\bf Proof.} By (\ref{3.7}) of Theorem 3.3 and Definition 4.1, we obtain
\bea i(c^{T+m_0+m})
&=& i(c^T)+(m+m_0)(i(c)+p_-+p_0-r ) + {{1-(-1)^{m_0+m}}\over 2}(q_0+q_+)  \nn\\
&&   +2\sum_{j=1}^r\left(E(\frac{(m_0+m+T)\th_j}{2\pi}) - E(\frac{T\th_j}{2\pi})\right)
       + 2\sum_{j=k_{\ast}+1}^{r_{\ast}}\vf(\frac{(m_0+m)\aa_j}{2\pi})  \nn\\
&=& i(c^T)+i(c^{m_0+m})+r+p_-+p_0+2(r_{\ast}-k_{\ast})+q_0+q_++2A-2k-(2A-2k)\nn\\
&&    +2\sum_{j=1}^r\left(E(\frac{(m_0+m+T)\th_j}{2\pi})
      -E(\frac{(m_0+m)\th_j}{2\pi})-E(\frac{T\th_j}{2\pi})\right)\nn\\
&=& i(c^T)+i(c^{m_0+m})+p(c)+2k-2A\nn\\
&&    +2\sum_{j=1}^k\left(E(\{\frac{m_0\th_j}{2\pi}\}
      +\{\frac{m\th_j}{2\pi}\}+\{\frac{T\th_j}{2\pi}\})
      -E(\{\frac{m_0\th_j}{2\pi}\}+\{\frac{m\th_j}{2\pi}\})
       -E(\{\frac{T\th_j}{2\pi}\})\right) \nn\\
&\ge& i(c^T)+i(c^{m_0+m})+p(c)-2\sum_{j=1}^kE(\{\frac{m_0\th_j}{2\pi}\}+\{\frac{m\th_j}{2\pi}\}) \nn\\
&\ge& i(c^T)+i(c^{m_0+m})+p(c)-4k\nn\\
&\ge& i(c^T)+p(c)+d,\qquad\forall \,m\ge 1.  \lb{4.26}\eea
This completes the proof of Lemma 4.4. \hfill\hb

\medskip

Next result generalizes Proposition 3.11 of \cite{LoD1} for rational closed geodesics to
irrational ones.

{\bf Theorem 4.5.} {\it For every orientable closed geodesic $c$ with $\hat{i}(c)>0$
on a Finsler manifold $(M,F)$ of dimension $d\ge 2$, denote the basic normal
form decomposition of the linearized Poincar\'e map $P_c$ of $c$ by (\ref{3.5}). Then
there exist an integer $A$ with $[(k+1)/2]\le A\le k$ and a subset $P$ of integers
$\{1, \ldots, k\}$ with $A$ integers such that for any small $\ep>0$
there exists a sufficiently large even integer $T\in n\N$ such that (\ref{4.5}) and (\ref{4.6})
and the following estimate hold}
\be  i(c^m)+\nu(c^m)\le i(c^T)+p(c)+d-3, \qquad\forall\,1\le m\le T-1. \lb{4.27}\ee

{\bf Proof.} When $k=0$, i.e., the closed geodesic $c$ is rational, this result was
proved in Proposition 3.11 of \cite{LoD1}, whose proof there in fact did not use
the fact $\;^{\#}\CG(M,F)=1$. Therefore here we only consider the case of $k\ge 1$.

On the one hand, by Theorem 3.3, for any $1\le m\le T-1$, we have
\bea
i(c^m)
&+& i(c^{T-m})+\nu(c^m) \nn\\
&=& i(c^{T}) + 2\sum_{j=1}^r\left(E(\frac{m\th_j}{2\pi})+E(\frac{(T-m)\th_j}{2\pi})
   -E(\frac{T\th_j}{2\pi})\right)- (r+p_-+p_0) \nn\\
&& - (-1)^m(q_0+q_+) +2\sum_{j=k_{\ast}+1}^{r_{\ast}}\vf(-\frac{m\aa_j}{2\pi})
   +2\sum_{j=k_{\ast}+1}^{r_{\ast}}\vf(\frac{m\aa_j}{2\pi})-2(r_{\ast}-k_{\ast}) \nn\\
&& +p_-+2p_0+p_++\frac{1+(-1)^m}{2}(q_-+2q_0+q_+)
  +2(r-k+r_{\ast}-k_{\ast}+r_0-k_0) \nn\\
&&-2[\sum_{j=k+1}^r\vf\left(\frac{m\th_j}{2\pi}\right)+
  \sum_{j=k_{\ast}+1}^{r_{\ast}}\vf\left(\frac{m\aa_j}{2\pi}\right)
    + \sum_{j=k_0+1}^{r_0}\vf\left(\frac{m\bb_j}{2\pi}\right)] \nn\\
&=& i(c^{T}) + 2\sum_{j=1}^r\left(E(\frac{m\th_j}{2\pi})+E(\frac{(T-m)\th_j}{2\pi})
 -E(\frac{T\th_j}{2\pi})\right)-r+p_0+p_+  \nn\\
&& - (-1)^m(q_0+q_+)+\frac{1+(-1)^m}{2}(q_-+2q_0+q_+)+2(r-k+r_0-k_0) \nn\\
&& +2\sum_{j=k_{\ast}+1}^{r_{\ast}}\vf(-\frac{m\aa_j}{2\pi})
  -2[\sum_{j=k+1}^r\vf\left(\frac{m\th_j}{2\pi}\right)+
   \sum_{j=k_0+1}^{r_0}\vf\left(\frac{m\bb_j}{2\pi}\right)],  \lb{4.28}\eea
where we have used the fact $T\in n\N$ is even and the fact $\nu(c)=p_-+p_++2p_0$
by the definitions of $p_{\ast}$s in Theorem 3.2.

Note that by (\ref{4.16}) we get
\bea
&& 2\sum_{j=1}^r\left(E(\frac{m\th_j}{2\pi}) + E(\frac{(T-m)\th_j}{2\pi})
    -E(\frac{T\th_j}{2\pi})\right)-2\sum_{j=k+1}^r\vf\left(\frac{m\th_j}{2\pi}\right)\nn\\
&&\qquad =2\sum_{j=1}^k\left(E(\frac{m\th_j}{2\pi})+E(\frac{(T-m)\th_j}{2\pi})
    -E(\frac{T\th_j}{2\pi})\right)\nn\\
&&\qquad\qquad+2\sum_{j=k+1}^r\left(E(\frac{m\th_j}{2\pi})
   +E(-\frac{m\th_j}{2\pi})\right)-2\sum_{j=k+1}^r\vf\left(\frac{m\th_j}{2\pi}\right)\nn\\
&&\qquad=2\sum_{j=1}^k\left(E(\{\frac{m\th_j}{2\pi}\})+E(\{\frac{T\th_j}{2\pi}\}
    -\{\frac{m\th_j}{2\pi}\})-E(\{\frac{T\th_j}{2\pi}\})\right)\nn\\
&&\qquad=2\sum_{j=1}^k\left(E(\{\frac{T\th_j}{2\pi}\}
    -\{\frac{m\th_j}{2\pi}\})\right)\nn\\
&&\qquad\le 2k.   \lb{4.29}\eea
Together with (\ref{4.28}), it yields
\bea i(c^m) &+& i(c^{T-m})+\nu(c^m) \nn\\
&&\le i(c^{T}) + 2k-r+p_0+p_++2(r-k+r_0-k_0)- (-1)^m(q_0+q_+)\nn\\
&& \quad+\frac{1+(-1)^m}{2}(q_-+2q_0+q_+)
   +2\sum_{j=k_{\ast}+1}^{r_{\ast}}\vf(\frac{m\aa_j}{2\pi})
   -2\sum_{j=k_0+1}^{r_0}\vf\left(\frac{m\bb_j}{2\pi}\right)\nn\\
&&=i(c^{T}) + r+p_0+p_++q_0+2(r_0-k_0)+\frac{1- (-1)^m}{2}q_+\nn\\
&&\qquad  \quad+\frac{1+(-1)^m}{2}q_-
   +2\sum_{j=k_{\ast}+1}^{r_{\ast}}\vf(\frac{m\aa_j}{2\pi})
   -2\sum_{j=k_0+1}^{r_0}\vf\left(\frac{m\bb_j}{2\pi}\right)\nn\\
&&\le i(c^T)+(p_-+p_0+q_0+q_++2r_{\ast}-2k_{\ast}+r+2A-2k)+p_++2(r_0-k_0)\nn\\
&&\qquad -p_--2(A-k)-\frac{1+(-1)^m}{2}q_++\frac{1+(-1)^m}{2}q_-
    -2\sum_{j=k_0+1}^{r_0}\vf\left(\frac{m\bb_j}{2\pi}\right)\nn\\
&&\le i(c^T)+p(c)+p_++q_-+2r_0-2(A-k)\nn\\
&&\le i(c^T)+p(c)+p_++q_-+2r_0+k, \qquad\forall\,1\le m\le T-1. \lb{4.30}\eea
In other words, we obtain
\be i(c^m)+\nu(c^m)\le i(c^T)+p(c)-i(c^{T-m})+p_++q_-+2r_0+k, \qquad 1\le m\le T-1. \lb{4.31}\ee

On the other hand, it follows from Theorem 3.8 that
\bea i(c^m)
&\le& i(c^T)-i(c)-p_0-p_-+r+q_0+q_++2(r_{\ast}-k_{\ast}) +2(A-k)  \nn\\
&=& i(c^T)+p(c)-i(c)-2(p_0+p_-),  \qquad \forall\,1\le m\le T-1.  \lb{4.32}\eea

Note that $p_+ + q_- + 2r_0+k\le d-1$ holds always in (\ref{4.31}) by (\ref{3.6})
with $d$ replaced by $d-1$. If
$p_+ + q_- + 2r_0+k \le d-3$, then (\ref{4.31}) yields (\ref{4.27}).
Therefore to continue our proof, it suffices to consider the following two distinct cases.

\medskip

{\bf Case 1.} {$p_+ + q_- + 2r_0+k = d-1$.}

\medskip

In this case, by (\ref{3.8}) and (\ref{3.10}) for all $m\ge 1$ we have
\be   \nu(c^m) \le \nu(c^n) = p_++q_-+2r_0 = d-1-k.   \lb{4.33}\ee
Thus together with (\ref{4.32}) it yields
\be i(c^m)+\nu(c^m)\le i(c^T)+p(c)-i(c)+d-1-k.  \lb{4.34}\ee
So in order to prove (\ref{4.27}), it suffices to consider the case of $i(c)=0$
and $k=1$. By the fact $i(c)=0$ and Proposition 3.4, we must have $q_-\in 2\N-1$
and thus $n\in 2\N$ by the definition of $n=n(c)$.

Therefore we have
\be  P_c\approx N_1(1,-1)^{\dm p_+}\dm N_1(-1,1)^{\dm q_-}
       \dm (\dm_{j=1}^{r_0}N_2(e^{\bb_j\sqrt{-1}},B_j))\dm R(\th_1),   \lb{4.35}\ee
where $\th_1/\pi\in (0,2)\bs\Q$. Thus by Theorem 3.3, we have
\bea
i(c^m) &=& -m + 2E(\frac{m\th_1}{2\pi}) - 1, \qquad \forall\; m\in\N, \lb{4.36}\\
\nu(c^n) &=& p_+ + q_- + 2r_0 = d-1-k=d-2.  \lb{4.37}\eea

When $m\in (\N\bs n\N)$, we have $\nu(c^m)< \nu(c^n)=\nu(c^T)$. Thus by
(\ref{4.32}) and (\ref{4.37}) we get
$$  i(c^m) + \nu(c^m) \le i(c^T) + p(c) + \nu(c^n) -1 = i(c^T) + p(c) + d-3, $$
i.e., (\ref{4.27}) holds.

Then for $1\le mn<T$ and the $T$ chosen above, by $\hat{i}(c)>0$ we get
\bea i(c^T)-i(c^{mn})&=& mn-T+2 \left(E(\frac{T\th_1}{2\pi})-E(\frac{mn\th_1}{2\pi})\right) \nn\\
&=& mn-T+(T-mn)\frac{\th_1}{\pi}+2 \left(\{\frac{mn\th_1}{2\pi}\}-\{\frac{T\th_1}{2\pi}\}\right) \nn\\
&=& (T-mn)\hat{i}(c)+2 \left(\{\frac{mn\th_1}{2\pi}\}-\{\frac{T\th_1}{2\pi}\}\right) \nn\\
&>& 2 \left(\{\frac{mn\th_1}{2\pi}\}-\{\frac{T\th_1}{2\pi}\}\right)   \nn\\
&>& -2.   \lb{4.38}\eea
Since both $n$ and $T$ are even, it follows from (\ref{4.36}) that $i(c^T)-i(c^{mn})$ is even.
Thus by the irrationality of $\frac{\th_1}{\pi}$ and (\ref{4.38}) we obtain
\be  i(c^T)\ge i(c^{mn}),\qquad\,\forall\, 1\le mn<T.  \lb{4.39}\ee
Then by the fact $p(c)=1$ and (\ref{4.37})-(\ref{4.39}) we have
\be  i(c^{mn}) + \nu(c^{mn}) \le  i(c^T)+\nu(c^n)
          \le i(c^T) + p(c)-1 + d-2=i(c^T)+p(c)+d-3.  \lb{4.40}\ee
That is, (\ref{4.27}) holds.

\medskip

{\bf Case 2.} $p_+ + q_- + 2r_0+k = d-2$

\medskip

In this case, $ p_- + p_0 + q_0 + q_+ + r-k + 2r_{\ast} + h_- + h_+=1$ by (\ref{3.6})
with $d$ replaced by $d-1$, which implies $r_{\ast}=0$. By Theorem 3.3 we have
\bea \nu(c^m)
&\le& \nu(c^n) \nn\\
&=& p_+ + q_- + 2r_0+(p_- +q_++ 2p_0 + 2q_0 + 2(r-k)) \nn\\
&=& d-2-k+(p_- +q_++ 2p_0 + 2q_0 + 2(r-k)),\quad\forall\,m\ge 1.   \lb{4.41}\eea

If $k\ge 3$, by (\ref{4.41}) it yields $\nu(c^m)\le d-3$, $\forall\,m\ge 1$. Thus
together with (\ref{4.32}), it yields (\ref{4.27}).

If $k=2$, by (\ref{4.32}) and (\ref{4.41}) it suffices to consider the following case
\be i(c)=p_0=p_-=q_+=h_+=h_-=r_{\ast}=0, \quad q_0+(r-k)=1, \quad k=2,  \lb{4.42}\ee
because otherwise (\ref{4.32}) would imply (\ref{4.27}) already.

Similarly, if $k=1$, by (\ref{4.31}), (\ref{4.32}) and (\ref{4.41}) it suffices to
consider the following case
\be i(c^{T-m})=p_0=p_-=h_+=h_-=r_{\ast}=0, \quad q_++q_0+(r-k)=1, \quad k=1. \lb{4.43}\ee
because otherwise (\ref{4.31}) and (\ref{4.32}) would imply (\ref{4.27}) already.

Now we consider (\ref{4.42}) and (\ref{4.43}) respectively.

\medskip

{\bf Case 2.1.} {\it (\ref{4.42}) happens.}

\medskip

Viewing $-I$ as $R(\pi)$ if $q_0=1$, it suffices to consider the case
$r-k=1$. Thus in addition to (\ref{4.42}) we have
\be  r=3, \qquad q_0=0.   \lb{4.44}\ee
Therefore we have
\be  P_c\approx N_1(1,-1)^{\dm p_+}\dm N_1(-1,1)^{\dm q_-}
    \dm (\dm_{j=1}^{r_0}N_2(e^{\bb_j\sqrt{-1}},B_j))\dm R(\th_1)\dm R(\th_2)\dm R(\th_3),   \lb{4.45}\ee
where $\th_1/\pi$ and $\th_2/\pi\in (0,2)\bs\Q$ and $\th_3/\pi\in(0,2)\cap\Q$. Thus by
Theorem 3.3, we have
\bea
i(c^m) &=& -3m + 2\sum_{j=1}^3E(\frac{m\th_j}{2\pi}) - 3, \qquad \forall\; m\in\N, \lb{4.46}\\
\nu(c^n) &=& p_+ + q_- + 2r_0 + 2 = d-2.  \lb{4.47}\eea
By the fact $i(c)=0$, (\ref{4.45}) and Proposition 3.4, there holds $q_-\in 2\N-1$.
By the definition of $n$, it further yields
\be n\in 2\N.  \lb{4.48}\ee

When $m\in (\N\bs n\N)$, we have $\nu(c^m)< \nu(c^n)=\nu(c^T)$. Thus by
(\ref{4.32}) and (\ref{4.47}) we get
$$  i(c^m) + \nu(c^m) \le i(c^T) + p(c) + \nu(c^n) -1 = i(c^T) + p(c) + d-3, $$
i.e., (\ref{4.27}) holds.

When $mn\in \N$ and $1\le mn<T$, then by (\ref{4.46}) and (\ref{4.48}) we have
$i(c^{T-mn})\in 2\N-1$. Therefore by (\ref{4.31}), for any $1\le mn<T$, we get
\bea  i(c^{mn}) + \nu(c^{mn})
&\le & i(c^T)+ p(c)-i(c^{T-mn})+d-2   \nn\\
&\le& i(c^T) + p(c) + d-3.  \lb{4.49}\eea
That is, (\ref{4.27}) holds.

\medskip

{\bf Case 2.2.} {\it (\ref{4.43}) happens.}

\medskip

Viewing $-I$ (or $N_1(-1,-1)$) as $R(\pi)$ if $q_0=1$ (or $q_+=1$),
although their nullity may be different by $1$ (cf. (\ref{4.53}) below), it suffices to
consider the case $r-k=1$. Thus in addition to (\ref{4.43}) we have
\be p(c)=2,\quad r=2, \quad q_+=q_0=0.   \lb{4.50}\ee
Therefore we have
\be  P_c\approx N_1(1,-1)^{\dm p_+}\dm N_1(-1,1)^{\dm q_-}
       \dm (\dm_{j=1}^{r_0}N_2(e^{\bb_j\sqrt{-1}},B_j))\dm R(\th_1)\dm R(\th_2),   \lb{4.51}\ee
where $\th_1/\pi\in (0,2)\bs\Q$ and $\th_2/\pi\in(0,2)\cap\Q$. Thus by Theorem 3.3, we have
\bea
i(c^m) &=& (i(c)-2)m + 2\sum_{j=1}^2E(\frac{m\th_j}{2\pi}) - 2, \qquad \forall\; m\in\N, \lb{4.52}\\
\nu(c^n) &=& p_+ + q_- + 2r_0 + 2(r-k) = d-1.  \lb{4.53}\eea

If $q_-\in 2\N-1$, by (\ref{4.51}) and Proposition 3.4, there holds
\be n\in 2\N \qquad \mbox{and}\qquad i(c^{T-m})\ge i(c)\in 2\N-1,\quad\forall\,1\le m\le T-1.\lb{4.54}\ee
Therefore, by (\ref{4.31}), (\ref{4.50}) and (\ref{4.53})-(\ref{4.54}) we get
\be  i(c^m)+\nu(c^m)\le i(c^T)+p(c)-i(c^{T-m})+d-2\le i(c^T)+p(c)+d-3,  \lb{4.55}\ee
That is, (\ref{4.27}) holds.

If $q_-\in2\N$, by (\ref{4.51}) and Proposition 3.4, it yields $i(c)\in 2\N_0$.
So it follows from (\ref{4.52}) that $i(c^m)\in 2\N_0$ for all $m\ge 1$. Let
$\frac{q}{p}=\frac{\th_2}{2\pi}$ with integers $p$ and $q$ satisfying $(p,q)=1$.

When $m\in(\N\bs p\N)$, we have $\nu(c^m)\le\nu(c^n)-2=d-3$ by (\ref{4.53}). Thus by
(\ref{4.32}) it yields
\be  i(c^m)+\nu(c^m)\le i(c^T)+p(c)+d-3.  \lb{4.56}\ee

When $m\in p\N$, then, similarly to (\ref{4.38}),
we can obtain $i(c^T)\ge i(c^m)$. Therefore, by (\ref{4.50}) and
(\ref{4.53}), we get
\be  i(c^m)+\nu(c^m)\le i(c^T)+\nu(c^n)\le i(c^T)+d-1+p(c)-2=i(c)+p(c)+d-3, \lb{4.57}\ee
That is, (\ref{4.27}) holds.

The proof is complete. \hfill\hb

\setcounter{equation}{0}
\section{Homological quasi-periodicity}

In this section, we study properties of homologies of energy level sets determined
by closed geodesics and establish certain periodicity of homological modules of
energy level set pairs when there exists only one prime closed geodesic.

For any $m\in\N$, denote the energy level $E(c^m)$ of $c^m$ by
\be \ka_m=E(c^m).  \lb{5.1}\ee
It is well known that $\ka_m=E(c^m)=m^2E(c)$ is strictly increasing to $+\infty$.
Set $\ka_0=0$. The next lemma follows from Theorem 3 of \cite{GrM1}, the Theorem
on p.367 of \cite{GrM2}, Lemma 3.1 to Theorem 3.7 of \cite{Lon4}, and
Theorem I.4.2 of \cite{Cha1}.

\medskip

{\bf Lemma 5.1.} (Lemma 4.2 of \cite{LoD1}) {\it Let $M=(M,F)$ be a Finsler
manifold with $\dim M<+\infty$. Let $c$ be a closed geodesic on $M$ each of whose
iteration $S^1\cdot c^m$ is an isolated critical orbit of $E$ in the loop space
$\Lm M$. Suppose that there are integers $m\in\N$ and $p\in 2\N_0$ such that
\be  i(c^m) = i(c) + p, \qquad \nu(c^m)=\nu(c).    \lb{5.2}\ee
Then the iteration map $\psi^m$ induces an isomorphism }
\be \psi^m_{\ast}: \ol{C}_{\ast}(E,c) \to \ol{C}_{\ast+p}(E,c^m).  \lb{5.3}\ee

\medskip

One of the key results in \cite{LoD1} is the homological isomorphism Theorem 4.3
there for rational closed geodesics. Below we redescribe this theorem and give more
details on two points for the proof given in \cite{LoD1}.

\medskip

{\bf Theorem 5.2.} (Theorem 4.3 of \cite{LoD1}) {\it Let $M=(M,F)$ be a
Finsler manifold possessing only one prime closed geodesic $c$ which is
rational and orientable. Let $n=n(c)$ be the analytical period of $c$.
Recall that by Theorem 3.7 of \cite{LoD1} there hold
\be  i(c^{m+n}) = i(c^m) + \ol{p}, \quad \nu(c^{m+n})=\nu(c^m), \qquad \forall\, m\in\N, \lb{5.4}\ee
where $\ol{p}=i(c^n)+p(c)$ is even. Then for any non-negative integers $b>a$ and any integer
$h\in\Z$, the iteration maps $\{\psi^m\}$ and inclusion maps of corresponding level sets
induce a map $f$ on singular chains which yields an isomorphism }
\be f_{\ast}: H_h(\ol{\Lm}^{\ka_b},\ol{\Lm}^{\ka_a})
   \to H_{h+\ol{p}}(\ol{\Lm}^{\ka_{n+b}},\ol{\Lm}^{\ka_{n+a}}).     \lb{5.5}\ee

{\bf Proof.} Proof of this theorem was given in \cite{LoD1} based on the above
Lemma 5.1. Note that an important condition in Lemma 5.1 is that the constant $p$
in (\ref{5.2}) should be even. In the applications of Lemma 5.1 (i.e., Lemma 4.2 of
\cite{LoD1}) in the proof of Theorem 4.3 in \cite{LoD1}, there are two points in its
Step 1 on which we did not give details on how to get this evenness condition. Below
we provide details of the proofs for these two points.

Because there is only one prime closed geodesic $c$ on $M$,
and $\ka_m=E(c^m)=m^2E(c)=m^2\ka_1>0$ is strictly increasing to $+\infty$, the critical
module of $E$ at $S^1\cdot c^m$ can be defined by
\be \ol{C}_j(E,c^m) = H_j(\ol{\Lm}^{\ka_m},\ol{\Lm}^{\ka_m\#})
       = H_j(\ol{\Lm}^{\ka_m},\ol{\Lm}^{\ka_{m-1}}),  \lb{5.6}\ee
where and below we denote by
\be   \ol{\Lm}^{\ka_m\#} \equiv \ol{\Lm^{\ka_m}\bs(S^1\cdot c^m)}.  \lb{5.7}\ee

Given a level set pair $(\ol{\Lm}^{\ka_p},\ol{\Lm}^{\ka_p\#})$ with $p\in\N$,
for any $\ga\in\Lm^{\ka_p}$ and $m\in\N$ we have
$$ E(\psi^m(\ga))=m^2E(\ga)\le m^2\ka_p=m^2E(c^p)=E(c^{mp})=\ka_{mp}. $$
Therefore the iteration map $\psi^m$ maps the level set $\ol{\Lm}^{\ka_p}$ into
$\ol{\Lm}^{\ka_{mp}}$. We denote the image of the pair $(\ol{\Lm}^{\ka_p},\ol{\Lm}^{\ka_p\#})$
under the iteration map $\psi^m$ by
\be (\ol{\Lm}^{\ka_p},\ol{\Lm}^{\ka_p\#})^m=(\psi^m(\ol{\Lm}^{\ka_p}),\psi^m(\ol{\Lm}^{\ka_p\#}))
        =(\ol{\psi^m(\Lm^{\ka_p})},\ol{\psi^m(\Lm^{\ka_p}\bs(S^1\cdot c^p))}\;). \lb{5.8}\ee

Note that we have (cf. (4.13)-(4.16) of \cite{LoD1})
\bea
&&  b=kn+q \qquad {\rm for\;some\;} k\in\N_0 \;\;{\rm and}\;\;0\le q\le n-1,  \lb{5.9}\\
&& i(c^b) = k\ol{p} + i(c^q),  \quad i(c^{n+b})=(k+1)\ol{p} + i(c^q),    \lb{5.10}\\
&& \nu(c^{n+b})=\nu(c^b)=\nu(c^q), \quad {\rm when}\;\;q\not= 0.  \lb{5.11}\eea

{\bf Point 1.} {\it The Proof of Case (i) with $\nu(c^b)=\nu(c)$ on Page 1787 of \cite{LoD1} }

Below (4.17) in Page 1787 of \cite{LoD1}, we have defined $\hat{p}=i(c^q)-i(c)$.

Now if $\hat{p}$ is even, then the constant $k\ol{p}+\hat{p}$ is even, and then we
can use Lemma 5.1 to get the isomorphism (4.23) in Page 1787 of \cite{LoD1}. Thus the
proof on Page 1787 for the Case (i) in \cite{LoD1} goes through.

Now if $\hat{p}$ is odd, then both $q$ and $n=n(c)$ must be even by (\ref{4.3}) of Lemma 4.2
and the definition of $n(c)$. Therefore $b$ is even by (\ref{5.9}). By (\ref{4.4}) of
Lemma 4.2 and (\ref{5.11}) we then obtain
$$  \nu(c^{n+b})=\nu(c^b)=\nu(c^q)\ge \nu(c^2)\ge \nu(c). $$
Thus equalities must hold here and we get
\be  \nu(c^{n+b})=\nu(c^b)=\nu(c^q) = \nu(c^2). \lb{5.12}\ee
We define $\td{p}=i(c^q)-i(c^2)$. Then $\td{p}$ is even by Lemma 4.2. We have also
\bea
i(c^b) &=& k\ol{p} + i(c^q) = k\ol{p} + \td{p} + i(c^2),  \lb{5.13}\\
i(c^{n+b}) &=& (k+1)\ol{p} + i(c^q) = (k+1)\ol{p} + \td{p} + i(c^2).  \lb{5.14}\eea

Thus we can replace (4.18)-(4.23) in \cite{LoD1} by the following arguments,
and obtain that the two iteration maps
\bea
\psi^{b/2}:&&(\ol{\Lm}^{\ka_2},\ol{\Lm}^{\ka_2\#})
          \to (\ol{\Lm}^{\ka_2},\ol{\Lm}^{\ka_2\#})^{b/2}
                  \subseteq (\ol{\Lm}^{\ka_b},\ol{\Lm}^{\ka_b\#}),   \lb{5.15}\\
\psi^{(n+b)/2}:&& (\ol{\Lm}^{\ka_2},\ol{\Lm}^{\ka_2\#}) \to (\ol{\Lm}^{\ka_2},\ol{\Lm}^{\ka_2\#})^{(n+b)/2}
                  \subseteq (\ol{\Lm}^{\ka_{n+b}},\ol{\Lm}^{\ka_{n+b}\#}), \lb{5.16}\eea
induce two isomorphisms on homological modules:
\bea
\psi^{b/2}_{\ast}:&&H_{h-k\ol{p}-\td{p}}(\ol{\Lm}^{\ka_2},\ol{\Lm}^{\ka_2\#})
          =\ol{C}_{h-k\ol{p}-\td{p}}(E,c^2) \to \ol{C}_{h}(E,c^b)
                 =H_h(\ol{\Lm}^{\ka_b},\ol{\Lm}^{\ka_b\#}),         \lb{5.17}\\
\psi^{(n+b)/2}_{\ast}:&& H_{h-k\ol{p}-\td{p}}(\ol{\Lm}^{\ka_2},\ol{\Lm}^{\ka_2\#})
          =\ol{C}_{h-k\ol{p}-\td{p}}(E,c^2)     \nn\\
     && \qquad\qquad \to \ol{C}_{h+\ol{p}}(E,c^{n+b})
      =H_{h+\ol{p}}(\ol{\Lm}^{\ka_{n+b}},\ol{\Lm}^{\ka_{n+b}\#}). \lb{5.18}\eea
Therefore the composed iteration map
\be f=\psi^{(n+b)/2}\circ \psi^{-b/2}: (\ol{\Lm}^{\ka_2},\ol{\Lm}^{\ka_2\#})^{b/2}
          \to (\ol{\Lm}^{\ka_2},\ol{\Lm}^{\ka_2\#})^{(n+b)/2}          \lb{5.19}\ee
is a homeomorphism and induces an isomorphism on homological modules:
\bea f_{\ast}:&& H_{h}(\ol{\Lm}^{\ka_b},\ol{\Lm}^{\ka_b\#})=\ol{C}_{h}(E,c^b) \to \ol{C}_{h+\ol{p}}(E,c^{n+b})
             =H_{h+\ol{p}}(\ol{\Lm}^{\ka_{n+b}},\ol{\Lm}^{\ka_{n+b}\#}), \lb{5.20}\eea
where we denote by $\psi^{-h}=(\psi^h)^{-1}$, the inverse map of $\psi^h$. Thus Theorem 5.2
holds in this case.

{\bf Point 2.} {\it The Proof of Case (iii-2) with $\nu(c^b)>\nu(c)$, $q>0$ in (\ref{5.9}),
and that there is some integer $t\in [1,q-1]$ such that $t|q$, $t|n$ and $\nu(c^t)=\nu(c^q)$ hold,
in Page 1789 of \cite{LoD1}. }

As in Page 1789 of \cite{LoD1}, let $s\in [1,q-1]$ be the minimal integer possessing the
property of the above integer $t$. Then $q=us$ and $n=vs$ hold for some $u, v\in\N$, and
as in \cite{LoD1} we obtain
\bea
&&  b = kn + (q-s) + s = (kv+u)s,        \lb{5.21}\\
&&  n+b = (k+1)n + (q-s) + s = ((k+1)v+u)s,   \lb{5.22}\\
&&  \nu(c^s) = \nu(c^q) = \nu(c^b) = \nu(c^{n+b}).  \lb{5.23}\eea
Let $\hat{p} = i(c^q) - i(c^s)$.

Now if $\hat{p}$ is even, then the constant $k\ol{p}+\hat{p}$ is even, and then we
can use Lemma 5.1 to get the isomorphism (4.46) in Page 1789 of \cite{LoD1}. Thus the
proof on Page 1789 for the Case (iii-2) in \cite{LoD1} goes through.

Now if $\hat{p}$ is odd, then $n=n(c)$ must be even by (\ref{4.3}) of Lemma 4.2 and
the definition of $n(c)$. By the same reason, $s$ and $q$ must have different parity.

Now if $s$ is even, then $q$ must be odd. Thus $b$ is odd by the evenness of $n$ and
(\ref{5.9}). This contradicts to (\ref{5.21}). Therefore $s$ must be odd and $q$ is even.

Because $s$ is odd, and both $s|q$ and $2|q$ hold, we have $(2s)|q$. Similarly $s|n$
and $2|n$ imply $(2s)|n$. Then by (\ref{5.9}) we obtain $(2s)|b$ and $(2s)|(n+b)$.

On the other hand, by (\ref{4.4}) of Lemma 4.2 and the fact $(2s)|q$, we obtain
$$  \nu(c^q)\ge \nu(c^{2s}) \ge \nu(c^s). $$
Together with (\ref{5.23}) we then obtain
\be  \nu(c^{n+b}) = \nu(c^b) = \nu(c^q) = \nu(c^{2s}) = \nu(c^s).  \lb{5.24}\ee
In this case we define $\td{p} = i(c^q)-i(c^{2s})$. Then $\td{p}$ is even by Lemma 4.2.
We have also
\bea
i(c^b) &=& k\ol{p} + i(c^q) = k\ol{p} + \td{p} + i(c^{2s}),  \lb{5.25}\\
i(c^{n+b}) &=& (k+1)\ol{p} + i(c^q) = (k+1)\ol{p} + \td{p} + i(c^{2s}).  \lb{5.26}\eea

Thus we can replace (4.41)-(4.46) in \cite{LoD1} by the following arguments, and
obtain from Lemma 5.1 that the two iteration maps
\bea
\psi^{b/(2s)}:&&(\ol{\Lm}^{\ka_{2s}},\ol{\Lm}^{\ka_{2s}\#})
          \to (\ol{\Lm}^{\ka_{2s}},\ol{\Lm}^{\ka_{2s}\#})^{b/(2s)}
                  \subseteq (\ol{\Lm}^{\ka_b},\ol{\Lm}^{\ka_b\#}),   \lb{5.27}\\
\psi^{(n+b)/(2s)}:&& (\ol{\Lm}^{\ka_{2s}},\ol{\Lm}^{\ka_{2s}\#})
                  \to (\ol{\Lm}^{\ka_{2s}},\ol{\Lm}^{\ka_{2s}\#})^{(n+b)/(2s)}
                  \subseteq (\ol{\Lm}^{\ka_{n+b}},\ol{\Lm}^{\ka_{n+b}\#}), \lb{5.28}\eea
induce two isomorphisms on homological modules:
\bea
\psi^{b/(2s)}_{\ast}:&&H_{h-k\ol{p}-\td{p}}(\ol{\Lm}^{\ka_{2s}},\ol{\Lm}^{\ka_{2s}\#})
          =\ol{C}_{h-k\ol{p}-\td{p}}(E,c^{2s})  \to \ol{C}_{h}(E,c^b)
                 =H_h(\ol{\Lm}^{\ka_b},\ol{\Lm}^{\ka_b\#}),         \lb{5.29}\\
\psi^{(n+b)/(2s)}_{\ast}:&& H_{h-k\ol{p}-\td{p}}(\ol{\Lm}^{\ka_{2s}},\ol{\Lm}^{\ka_{2s}\#})
          =\ol{C}_{h-k\ol{p}-\td{p}}(E,c^{2s})    \nn\\
  && \qquad\qquad \to \ol{C}_{h+\ol{p}}(E,c^{n+b})
      =H_{h+\ol{p}}(\ol{\Lm}^{\ka_{n+b}},\ol{\Lm}^{\ka_{n+b}\#}). \lb{5.30}\eea
Therefore the composed iteration map
\be f=\psi^{(n+b)/(2s)}\circ \psi^{-b/(2s)}: (\ol{\Lm}^{\ka_{2s}},\ol{\Lm}^{\ka_{2s}\#})^{b/(2s)}
          \to  (\ol{\Lm}^{\ka_{2s}},\ol{\Lm}^{\ka_{2s}\#})^{(n+b)/(2s)}     \lb{5.31}\ee
is a homeomorphism and induces an isomorphism on homological modules:
\be
f_{\ast}: H_{h}(\ol{\Lm}^{\ka_b},\ol{\Lm}^{\ka_b\#})=\ol{C}_{h}(E,c^b)  \to \ol{C}_{h+\ol{p}}(E,c^{n+b})
      =H_{h+\ol{p}}(\ol{\Lm}^{\ka_{n+b}},\ol{\Lm}^{\ka_{n+b}\#}).  \lb{5.32}\ee
Thus Theorem 5.2 holds in this case too.

Now the rest part of the proof of Theorem 4.3 of \cite{LoD1} yields Theorem 5.2. \hfill\hb

\medskip

The above homological isomorphism theorem is for rational closed geodesics. Our next
result generalizes it to irrational closed geodesics, and will play a crucial role in the
proofs of Theorems 1.1 and 1.2. Here the quasi-periodicity which we established in the
above Theorem 4.3 is crucial in the proof.

\medskip

{\bf Theorem 5.3.} {\it Let $(M,F)$ be a Finsler manifold possessing only one
prime closed geodesic $c$ which is orientable and $n=n(c)$ be the analytical
period of $c$. Recall that by Theorem 4.3 there exists an even integer $T\in n\N$ such
that for $m_0=m_0(c)$ given by Definition 4.1 there hold
\be  i(c^{m+T}) = i(c^m) + \ol{p}, \quad \nu(c^{m+T})=\nu(c^m), \qquad \forall\, 1\le m\le m_0,
       \lb{5.33}\ee
where $\ol{p}=i(c^T)+p(c)$. Then we can further require $T\in (m_0!n)\N$ such that for any
non-negative integers $a$ and $b$ satisfying $0<a<b\le m_0$, the iteration maps $\{\psi^m\}$
and inclusion maps of corresponding level sets induce a map $f$ on singular chains which
yields an isomorphism }
\be f_{\ast}: H_h(\ol{\Lm}^{\ka_b},\ol{\Lm}^{\ka_a})
   \to H_{h+\ol{p}}(\ol{\Lm}^{\ka_{T+b}},\ol{\Lm}^{\ka_{T+a}}),\qquad \forall\, h\in \Z.   \lb{5.34}\ee

{\bf Proof.} Here we follow the main ideas from pages 1786-1792 of \cite{LoD1}.

{\bf Step 1.} {\it The isomorphism in the case of $\;b-a=1$. }

In this case $\ka_a$ and $\ka_b$ are the only two critical values in $[\ka_a, \ka_b]$
and so are $\ka_{T+a}$ and $\ka_{T+b}$ in $[\ka_{T+a}, \ka_{T+b}]$. Then, note that
$0\le a<b\le m_0$, by (\ref{5.33}) it yields
\be i(c^{T+b}) = \ol{p} + i(c^b),\qquad \nu(c^{T+b})=\nu(c^b). \lb{5.35}\ee
Here we require the even integer $T$ chosen by Theorems 4.3-4.5 to further satisfy
$T\in (m_0!n)\N$. Thus it yields
\be   b\,|(T+b), \qquad\forall\, 0\le a<b\le m_0.  \lb{5.36}\ee

Note first that $\ol{p}=i(c^T)+p(c)$ is always even by (B) of Theorem 4.3.
Therefore by (\ref{5.35}) we get
\be  \ep(c^{T+b})=(-1)^{i(c^{T+b})-i(c)}=(-1)^{i(c^b)-i(c)}=\ep(c^b).  \lb{5.37}\ee

For any $0\le a<b\le m_0$, since $\nu(c^{T+b})=\nu(c^b)$ holds in (\ref{5.35})
and $b|(T+b)$ holds in (\ref{5.36}), it follows from Lemma 2.1 and (iii) of Lemma 2.2 that
\bea H_{h}(\ol{\Lm}^{\ka_b},\ol{\Lm}^{\ka_b\#})
&=& \ol{C}_{h}(E,c^b)\nn\\
&=& H_{h-i(c^b)}(N_{c^b}^-\cup\{c^b\},N_{c^b}^-)^{\ep(c^b)\Z_b}\nn\\
&=& H_{h-i(c^b)}(N_{c^{T+b}}^-\cup\{c^{T+b}\},N_{c^{T+b}}^-)^{\ep(c^{T+b})\Z_{T+b}}\nn\\
&=& H_{h+\ol{p}-i(c^{T+b})}(N_{c^{T+b}}^-\cup\{c^{T+b}\},N_{c^{T+b}}^-)^{\ep(c^{T+b})\Z_{T+b}}\nn\\
&=& \ol{C}_{h+\ol{p}}(E,c^{T+b})\nn\\
&=& H_{h+\ol{p}}(\ol{\Lm}^{\ka_{T+b}},\ol{\Lm}^{\ka_{T+b}\#}). \lb{5.38}\eea
Here we used Lemma 2.1 in the second and fifth equalities, (iii) of Lemma 2.2
and (\ref{5.35})-(\ref{5.36}) in the third one and (\ref{5.35}) in the fourth one.

The case of $b-a=1$ is proved.

\medskip

{\bf Step 2.} {\it The induction argument for general $b>a$.}

\medskip

Now we can follow precisely the proof in the Step 2 on pages 1789-1792 of Theorem 4.3 of
\cite{LoD1} and complete the proof of Theorem 5.3 here. Thus we omit all these details
here. \hfill\hb

\medskip

Next we generalize the Proposition 5.1 of \cite{LoD1} for rational closed geodesics to
irrational ones. Here we denote by $\Q^m$ the $m$ times of the module instead of using
the notation $m\Q$ in order to make the text clearer.

\medskip

{\bf Theorem 5.4.} {\it Let $c$ be the only one prime closed geodesic on a compact Finsler
manifold $(M,F)$ of dimension $d$. Suppose $c$ is orientable. Let $n=n(c)$ be the
analytical period of $c$ and $m_0=m_0(c)$ be given by Definition 4.1. Let
$T\in (m_0!n)\N$ be the even integer given by Theorems 4.3-4.5 and 5.3. Denote by
$X_j=H_j(\ol{\Lm},\ol{\Lm}^T)=\Q^{x_j}$ for all $j\in\Z$. Then there holds
\be  x_j = b_{j-i(c^T)-p(c)} \qquad \forall\;0\le j\le i(c^T)+p(c)+d-2, \lb{5.39}\ee
where $b_j$'s are the Betti numbers of the loop space $\Lm M$ defined in Section 2. }

\medskip

{\bf Proof.} Let $R=i(c^T)$. Firstly we fix an integer $j\le R+p(c)+d-2$. Because there is
only one prime closed geodesic $c$ on $M$, there holds $\hat{i}(c)>0$. Thus we have
$i(c^m)\to +\infty$ as $m\to +\infty$. According to Definition 4.1 and Lemma 4.4 we have
$$  i(c^m)\ge R+p(c)+d,  \qquad \forall\;m\ge T+m_0.  $$
It then implies
$$  \ol{C}_q(E,c^m) = 0, \qquad \forall\; m\ge T+m_0, \quad q\le j+1=R+p(c)+d-1.  $$
Therefore by Theorem II.1.5 on page 89 of \cite{Cha1} we obtain
\be  H_q(\ol{\Lm},\ol{\Lm}^{\ka_m})=0 \qquad \forall\;m\ge T+m_0,\;\; 0\le q\le j+1. \lb{5.40}\ee
Thus the exact sequence of the triple $(\ol{\Lm},\ol{\Lm}^{\ka_{m_0+T}},\ol{\Lm}^{\ka_T})$ yields
\be  0=H_{j+1}(\ol{\Lm},\ol{\Lm}^{\ka_{m_0+T}})\to H_{j}(\ol{\Lm}^{\ka_{m_0+T}},\ol{\Lm}^{\ka_T})
      \to H_{j}(\ol{\Lm},\ol{\Lm}^{\ka_T})\to H_{j}(\ol{\Lm},\ol{\Lm}^{\ka_{m_0+T}})=0,  \lb{5.41}\ee
which then implies the isomorphism:
\be H_j(\ol{\Lm}^{\ka_{m_0+T}},\ol{\Lm}^{\ka_T}) = H_j(\ol{\Lm},\ol{\Lm}^{\ka_T})=\Q^{x_j}. \lb{5.42}\ee

On the other hand, by Theorem 5.3 we obtain an isomorphism:
\be H_{j-R-p(c)}(\ol{\Lm}^{\ka_{m_0}},\ol{\Lm}^0) = H_j(\ol{\Lm}^{\ka_{m_0+T}},\ol{\Lm}^{\ka_T}),
           \quad\forall\,j\le R+p(c)+d-2.  \lb{5.43}\ee

Fix an integer $l\le d-2$. By Definition 4.1 we have $i(c^m)\ge d+4k$ for all $m\ge m_0$,
which implies
\be  H_q(\ol{\Lm},\ol{\Lm}^{\ka_m})=0 \qquad \forall\;m\ge m_0,\;\; 0\le q\le l+1. \lb{5.44}\ee

Then the exact sequence of the triple $(\ol{\Lm},\ol{\Lm}^{\ka_{m_0}},\ol{\Lm}^0)$ yields
\bea  0
&=& H_{j-R-p(c)+1}(\ol{\Lm},\ol{\Lm}^{\ka_{m_0}})\to H_{j-R-p(c)}(\ol{\Lm}^{\ka_{m_0}},\ol{\Lm}^0)   \nn\\
&& \qquad\qquad  \to H_{j-R-p(c)}(\ol{\Lm},\ol{\Lm}^0)\to H_{j-R-p(c)}(\ol{\Lm},\ol{\Lm}^{\ka_{m_0}})=0,
\quad\forall\,j\le R+p(c)+d-2.               \nn\eea
It then implies the isomorphism:
\be  H_{j-R-p(c)}(\ol{\Lm}^{\ka_{m_0}},\ol{\Lm}^0) = H_{j-R-p(c)}(\ol{\Lm},\ol{\Lm}^0)
 = \Q^{b_{j-R-p(c)}},\quad\forall\,j\le R+p(c)+d-2.          \lb{5.45}\ee
Therefore (\ref{5.42})-(\ref{5.43}) and (\ref{5.45}) yield the claim (\ref{5.39}). \hfill\hb

\medskip

Next we generalize the Theorems 5.2 of \cite{LoD1} for the rational closed geodesics to
irrational ones.

{\bf Theorem 5.5.} {\it Let $(M,F)$ be a compact simply connected $dh$-dimensional Finsler
manifold with $H^{\ast}(M,\Q)=T_{d,h+1}(x)$ for some integers $d\ge 2$ and $h\ge 1$.
Suppose $c$ is the only one prime closed geodesic on $M$, and let $\mu=p(c)+dh-3$. Denote
by $n=n(c)$ and $m_0=m_0(c)$ given by (\ref{4.1}) and Definition 4.1. Then there exist an
even integer $T\in (m_0!n)\N$ and an integer $\ka\ge 0$ such that
\be B(d,h)(i(c^T) + p(c)) + (-1)^{\mu+i(c^T)}\ka = \sum_{j=\mu-p(c)+1}^{i(c^T)+\mu}(-1)^j b_j,
                 \lb{5.46}\ee
where $B(d,h)$ is given in Lemma 2.4. }

\medskip

{\bf Proof.} Note first that by Proposition 3.5 and Remark 3.6 the closed geodesic $c$
on $M$ is orientable because $M$ is simply connected.

Since there exists only one prime closed geodesic $c$, it follows that
$\hat{i}(c)>0$ and $0\le i(c)\le d-1$. Specially by Lemma 2.4 we obtain
\be  \hat{i}(c)\in \Q.  \lb{5.47}\ee
Let
\be d_j = k_j^{\ep(c^n)}(c^n), \qquad \forall j\in\Z.   \lb{5.48}\ee
Then by the definition of $n=n(c)$, Lemma 2.3 and (\ref{4.2}) we obtain
\be k_j^{\ep(c^{mn})}(c^{mn}) = d_j, \qquad \forall j\in\Z, \quad m\in\N.   \lb{5.49}\ee

Fix $T\in (m_0!n)\N$ to be an even integer determined by Theorems 4.3-4.5, 5.3, and 5.5.
Specially we require that this $T$ makes (\ref{4.12}) hold.

Then we claim the following four conditions hold:
\bea
&&  i(c^{m+T}) = i(c^T) + i(c^m) + p(c), \qquad \forall\; 1\le m\le m_0,  \lb{5.50}\\
&&  i(c^m)+\nu(c^m) \le i(c^T)+\mu,  \qquad \forall\; 1\le m<T, \lb{5.51}\\
&&  d_j = 0, \qquad \forall\;j\ge \mu+2,  \lb{5.52}\\
&&  H_{i(c^T)+\mu+1}(\ol{\Lm},\ol{\Lm}^{\ka_T})=0.  \lb{5.53}\eea

In fact, (\ref{5.50}) follows from (A) of Theorem 4.3, and (\ref{5.51}) follows from
Theorem 4.5.

Note that if $k\ge 1$ in Theorem 4.3, there holds $A\ge 1$ by Proposition 3.7. Thus for
$j\ge\mu+2=p(c)+dh-1$, it yields $j>\nu(c^n)$ by (C) of Theorem 4.3, which implies
that (\ref{5.52}) holds. If $k=0$, then (\ref{5.52}) was proved in the proof of Theorem
6.1 of \cite{LoD1} when verifying the condition (5.11) there via Hingston's Theorem
of \cite{Hin2} (cf. Theorem 4.1 of \cite{LoD1}).

Note that $i(c^T)+\mu+1=i(c^T)+p(c)+dh-2$ holds. So by Theorem 5.4, Lemmas 2.5 and 2.6, we
obtain
$$  H_{i(c^T)+\mu+1}(\ol{\Lm},\ol{\Lm}^{\ka_T})=b_{dh-2}=0.  $$
Thus (\ref{5.53}) holds, and the proof of the four conditions (\ref{5.50})-(\ref{5.53})
is complete.

Let $R=i(c^T)$. Then by (\ref{5.50})-(\ref{5.53}) and Lemma 4.4 we obtain the
following distribution diagram (\ref{5.54}) of $\dim\ol{C}_j(E,c^m)$ for any
$j\ge 0$ and $m\ge 1$.
{\footnotesize\bea
\begin{tabular}{c|ccc ccc ccc ccc ccc}
$\cdots$& & & & & &&&&&&$\ast$&$\cdots$\\
$T+m_0+1$& & & & & &&&&&&$\ast$&$\cdots$\\
$T+m_0$& & & & & &$\ast$&$\cdots$&$\cdots$&$\cdots$&$\cdots$&$\ast$&$\cdots$\\
$\cdots$& & & & & &$\cdots$&$\cdots$&$\cdots$&$\cdots$&$\cdots$&$\cdots$&$\cdots$\\
$T+1$& & & & & &$\ast$&$\cdots$&$\cdots$&$\cdots$&$\cdots$&$\ast$&$\cdots$\\
$T$& & &$0$&$d_0$&$\cdots$&$d_{p(c)}$&$\cdots$&$d_{\mu}$&$d_{\mu+1}$&$d_{\mu+2}$&$0$&$0$\\
$T-1$&$\ast$&$\cdots$&$\cdots$&$\cdots$&$\cdots$&$\cdots$&$\cdots$&$\ast$  \\
$\cdots$&$\cdots$&$\cdots$&$\cdots$&$\cdots$&$\cdots$&$\cdots$&$\cdots$&$\cdots$ \\
$1$&$\ast$&$\cdots$&$\cdots$&$\cdots$&$\cdots$&$\cdots$&$\cdots$&$\ast$ \\
\cline{1-13}
$m$ in $c^m$&$c_{0}$&$\cdots$&$c_{R-1}$&$c_{R}$
       &$\cdots$&$c_{R+p(c)}$&$\cdots$&$c_{R+\mu}$&$c_{R+\mu+1}$
                  &$c_{R+\mu+2}$&$c_{R+\mu+3}$&$\cdots$ \\
\end{tabular}  \lb{5.54}\eea}

Here as coordinates of the diagram (\ref{5.54}), the first left column lists the
iteration time $m$ of $c^m$ starting from $1$ to $T+m_0+1$ and upwards, and the
first row from below lists the dimensions $c_j=\dim\ol{C}_j(E,c^{\ast})$ of the
$S^1$-equivariant critical module $\ol{C}_j=\ol{C}_j(E,c^{\ast})$ from $j=0$ to
$j=R+\mu+3$ and rightwards. The entry $D_j(c^m)$ in this diagram at $m$-th row and
$j$-th column is given by $D_j(c^m)=\dim \ol{C}_j(E,c^m)$. Here $d_j=\dim \ol{C}_j(E,c^T)$s
are shown in this diagram. $\ast$s and Dots in this diagram indicate entries which may not
be zero whose precise values depend on $\dim \ol{C}_j(E,c^m)$. Entries on the empty places
in the diagram are all $0$.

Now the proof is similar to that of Theorem 5.2 in \cite{LoD1} (cf. pages 1795-1799 for
more details). Here for reader's conveniences, we include certain details of the proof
here.

Denote by $\ka_m=E(c^m)$ for $m\ge 1$. As in the Step 1 of the proof of Theorem 5.2 of
\cite{LoD1}, for $j\in\Z$, we denote by
\be U_j=H_j(\ol{\Lm}^{\ka_T},\ol{\Lm}^0)=\Q^{u_j}, \quad B_j=H_j(\ol{\Lm},\ol{\Lm}^0)=\Q^{b_j},
          \quad X_j=H_j(\ol{\Lm},\ol{\Lm}^{\ka_T})=\Q^{x_j}.  \lb{5.55}\ee

Then the long exact sequence of the triple
$(\ol{\Lm},\ol{\Lm}^{\ka_T},\ol{\Lm}^0)$ yields the following diagram:
\bea
\begin{tabular}{ccc ccc ccc ccc ccc}
$X_{R+\mu+1}$&$\to$&$U_{R+\mu}$&$\to$&$B_{R+\mu}$&$\to$&$X_{R+\mu}$
                     &$\to$&$\cdots$&$\to$&$U_0$&$\to$&$B_0$&$\to$&$X_0$ \\
$\parallel$&   &$\parallel$&   &$\parallel$&   &$\parallel$&  & & &$\parallel$& &$\parallel$& &$\parallel$ \\
$0$& &$\Q^{u_{R+\mu}}$&  &$\Q^{b_{R+\mu}}$& &$\Q^{x_{R+\mu}}$& &$\cdots$& &$\Q^{u_0}$& &$0$& &$\;0$, \\
\end{tabular}  \nn\eea
where $X_{R+\mu+1}=0=X_0$ follows from (\ref{5.53}), Theorem 5.4 and Lemmas 2.5 and 2.6.
$B_0=0$ follows from Lemmas 2.5 and 2.6. Then this long exact sequence yields
\be  0 = \sum_{j=0}^{R+\mu}(-1)^j(u_j - b_j + x_j).  \lb{5.56}\ee

Because $T\ge 2$, for $j\in\Z$ besides $U_j$ defined in (\ref{5.55}) we denote by
$$ V_j=H_j(\ol{\Lm}^{\ka_{T-1}},\ol{\Lm}^0)=\Q^{v_j},
           \quad E_j=H_j(\ol{\Lm}^{\ka_T},\ol{\Lm}^{\ka_{T-1}})=\Q^{e_j}. $$
Then the exact sequence of the triple $(\ol{\Lm}^{\ka_T},\ol{\Lm}^{\ka_{T-1}},\ol{\Lm}^0)$ and
the diagram (\ref{5.54}) yield the following diagram:
\bea
\begin{tabular}{ccc ccc ccc ccc ccc}
$V_{R+\mu+1}$&$\to$&$U_{R+\mu+1}$&$\to$&$E_{R+\mu+1}$&$\to$&$V_{R+\mu}$&$\to$&$\cdots$ \\
$\parallel$&   &$\parallel$&   &$\parallel$&   &$\parallel$&  & &      \\
$0$& &$\Q^{u_{R+\mu+1}}$&  &$\Q^{e_{R+\mu+1}}$& &$\Q^{v_{R+\mu}}$& &$\cdots$&  \\
&&&&&&&& \\
 &$\to$&$V_{R}$&$\to$&$U_{R}$&$\to$&$E_{R}$&$\to$&$\cdots$ \\
 &  &$\parallel$&   &$\parallel$&   &$\parallel$&  &     \\
 & &$\Q^{v_{R}}$&  &$\Q^{u_{R}}$& &$\Q^{e_{R}}$& &$\cdots$ \\
&&&&&&&& \\
 &$\to$&$V_{0}$&$\to$&$U_{0}$&$\to$&$E_{0}$&$\to$&$0$  \\
 &  &$\parallel$&   &$\parallel$&   &$\parallel$&      \\
 & &$\Q^{v_{0}}$&  &$\Q^{u_{0}}$& &$\Q^{e_0}$  \\
\end{tabular}  \nn\eea
where $V_{R+\mu+1}=0$ follows from (\ref{5.51}) and the diagram (\ref{5.54}). Then this long exact
sequence yields
\be \sum_{j=0}^{R+\mu}(-1)^ju_j = (-1)^{R+\mu}u_{R+\mu+1} + \sum_{j=0}^{R+\mu}(-1)^jv_j
           + \sum_{j=0}^{R+\mu+1}(-1)^je_j.         \lb{5.57}\ee
Note that by (\ref{5.49}) we have
\be  e_j = \left\{\matrix{
             d_{j-R}, & \quad {\rm for}\;\;R \le j\le R+\mu+1, \cr
             0, & \quad {\rm otherwise}.  \cr}\right.  \lb{5.58}\ee
Thus we obtain
\be \sum_{j=0}^{R+\mu}(-1)^ju_j = (-1)^{R+\mu}u_{R+\mu+1} + \sum_{j=0}^{R+\mu}(-1)^jv_j
           + \sum_{j=0}^{\mu+1}(-1)^{R+j}d_j.     \lb{5.59}\ee

Now combining (\ref{5.56}) and (\ref{5.59}) we obtain
\be  0 = \sum_{j=0}^{R+\mu}(-1)^jv_j + \sum_{j=0}^{\mu+1}(-1)^{R+j}d_j
      - \sum_{j=0}^{R+\mu}(-1)^jb_j
      +\sum_{j=0}^{R+\mu}(-1)^jx_j + (-1)^{R+\mu}u_{R+\mu+1}.     \lb{5.60}\ee

Now as in \cite{LoD1}, we can apply the procedure above to decrease the level sets
one by one by induction. In this way, each time we pass through a critical level
$E(c^m)$ with $m\le T$, the term $\sum_{j=0}^{R+\mu}(-1)^jv_j$ on the right hand
side of (\ref{5.60}) will be replaced by the sum of a similar alternating sum of
dimensions of homological modules of a new lower level set pair
$(\ol{\Lm}^{\ka_{m-1}},\ol{\Lm}^0)$ and a term $\sum_{j=0}^{\nu(c^m)}(-1)^{i(c^m)+j}k_j^{\ep(c^m)}(c^m)$.
Here the sign of $i(c^m)$ indicates the parity of the number of column in which the term
$k_0^{\ep(c^m)}(c^m)$ appears. Then by induction from (\ref{5.56})-(\ref{5.60}) repeating
the proof of Theorem 5.2 in \cite{LoD1} by using the above diagram (\ref{5.54}) and our
Theorem 5.4, similarly to (5.22) of \cite{LoD1} we obtain
\bea  0
&=& \frac{T}{n}\sum_{j=0}^{\mu+1}(-1)^{i(c^n)+j}d_j
      + \frac{T}{n}\sum_{m=1}^{n-1}\sum_{j=0}^{\nu(c^m)}(-1)^{i(c^m)+j}k_j^{\ep(c^m)}(c^m)   \nn\\
&& \qquad - \sum_{j=0}^{R+\mu}(-1)^j b_j + \sum_{j=0}^{R+\mu}(-1)^jb_{j-R-p(c)}
      + (-1)^{R+\mu}u_{R+\mu+1}.   \lb{5.61}\eea
Note that in the proof of (\ref{5.61}), the facts that $T$ is an integer multiple of
$n(c)$, the $n(c)$-periodicity of critical modules in iterates given by Lemma 2.3,
(\ref{4.2}) and (\ref{5.49}) are crucial.

Now similarly to (5.23) of \cite{LoD1} we can apply the mean index identity
Lemma 2.4 to further obtain
\bea B(d,h)n\hat{i}(c)
&=& \sum_{1\le m\le n\atop 0\le j\le 2dh-2}(-1)^{i(c^m)+j}k_j^{\epsilon_j}(c^m)  \nn\\
&=& \sum_{1\le m\le n-1\atop 0\le j\le 2dh-2}(-1)^{i(c^m)+j}k_j^{\epsilon_j}(c^m)
   +\sum_{j=0}^{\nu(c^n)}(-1)^{i(c^n)+j}k_j^{\epsilon_n}(c^n)  \nn\\
&=& \sum_{m=1}^{n-1}\sum_{j=0}^{\nu(c^m)}(-1)^{i(c^m)+j}k_j^{\epsilon_j}(c^m)
   +\sum_{j=0}^{\nu(c^n)}(-1)^{i(c^n)+j}k_j^{\epsilon_n}(c^n)  \nn\\
&=& \sum_{m=1}^{n-1}\sum_{j=0}^{\nu(c^m)}(-1)^{i(c^m)+j}k_j^{\epsilon_j}(c^m)  \nn\\
& & \qquad +\sum_{j=0}^{\mu+1}(-1)^{i(c^n)+j}d_j+\sum_{j=\mu+2}^{\nu(c^n)}(-1)^{i(c^n)+j}d_j  \nn\\
&=& \sum_{m=1}^{n-1}\sum_{j=0}^{\nu(c^m)}(-1)^{i(c^m)+j}k_j^{\epsilon_j}(c^m)
   +\sum_{j=0}^{\mu+1}(-1)^{i(c^n)+j}d_j,  \lb{5.62}\eea
where we have used the condition $d_j=0$ for all $j\ge\mu+2$ of (5.52) in the last equality.

Now by (D) of Theorem 4.3, the rationality (\ref{5.47}) of $\hat{i}(c)$, (\ref{5.61})
and (\ref{5.62}) we obtain
\bea  0
&=& B(d,h)T\hat{i}(c) - \sum_{j=0}^{R+\mu}(-1)^jb_j
      + \sum_{j=0}^{R+\mu}(-1)^jb_{j-R-p(c)} + (-1)^{R+\mu}u_{R+\mu+1}   \nn\\
&=& B(d,h)(R + p(c)) - \sum_{j=0}^{R+\mu}(-1)^jb_j
      + \sum_{j=0}^{R+\mu}(-1)^jb_{j-R-p(c)} + (-1)^{R+\mu}u_{R+\mu+1}   \nn\\
&=& B(d,h)(R + p(c)) - \sum_{j=\mu-p(c)+1}^{R+\mu}(-1)^jb_j
        + (-1)^{R+\mu}u_{R+\mu+1}.         \lb{5.63}\eea
That is, (\ref{5.46}) holds with $\ka=u_{R+\mu+1}\ge 0$.

This completes the proof of Theorem 5.5. \hfill\hb

\setcounter{equation}{0}
\section{Proofs of Theorems 1.1 and 1.2}

In this section, we will follow ideas from Section 6 of \cite{LoD1} and Section 4 of
\cite{DuL3} to give the proofs of Theorems 1.1 and 1.2 via replacing $n=n(c)$ by the
integer $T$ obtained by Theorems 4.3, 4.5, 5.3 and 5.5, and modifying related arguments
using our above results. For reader's conveniences and completeness, we give all the
details here.

\medskip

{\bf Proof of Theorem 1.1.} Let $M$ be a compact simply connected manifold of dimension
not less than $2$ with a Finsler metric $F$. By Theorems A and B in the Section 1, it
suffices to assume that the condition (\ref{1.3}) on $M$ holds, i.e.,
$$  H^*(M;\Q)\cong T_{d,h+1}(x)=\Q[x]/(x^{h+1}=0)   $$
with a generator $x$ of degree $d\ge 2$ and hight $h+1\ge 2$.

We prove the theorem by contradiction. Thus we assume that there exists only one prime
closed geodesic $c$ on $(M,F)$. To generate the non-trivial $H_{d-1}(\Lm M/S^1,\Lm M^0/S^1;\Q)$
(cf. Lemmas 2.5 and 2.6), this $c$ must satisfy
\be  0\le i(c)\le d-1, \quad \hat{i}(c)>0, \quad \hat{i}(c)\in \Q.  \lb{6.1}\ee
where the last conclusion follows from Lemma 2.4.

For the analytic period $n=n(c)$ and $m_0=m_0(c)$ given by Definition 4.1, fix a large
even integer $T\in (m_0!n)\N$ determined by Theorems 4.3, 4.5, 5.3 and 5.5. Then by (\ref{6.1})
and (D) of Theorem 4.3 we have
$$  i(c^T) + p(c) = T \hat{i}(c)>0.  $$
Note that $i(c^T)=p(c)$ $(\mod\;2)$ by (B) of Theorem 4.3, so we obtain
\be  i(c^T) + p(c) \in 2\N.   \lb{6.2}\ee
Let $\mu=p(c)+(dh-3)$. Then by (\ref{6.2}) we have
\be  i(c^T) + \mu \ge dh-1 \ge 1, \qquad i(c^T) + \mu \in 2\N_0+(dh-1).  \lb{6.3}\ee
Then by Theorem 5.5, we obtain for some integer $\ka\ge 0$:
\be B(d,h)(i(c^T) + p(c))  + (-1)^{i(c^T)+\mu}\ka
          = \sum_{j=\mu-p(c)+1}^{i(c^T)+\mu}(-1)^jb_j.  \lb{6.4}\ee

Note that when $d$ is odd, then $h=1$ by Remark 2.5 of \cite{Rad1}. And when $h=1$, $M$
is rationally homotopic to the sphere $S^d$. So we can classify the manifolds $M$ satisfying
(\ref{1.3}) into two classes according to the parity of $d$, and continue our proof
correspondingly.

\medskip

{\bf Case 1.} {$d\ge 2$ is even and $h\ge 1$.}

\medskip

Note that, in this case, $i(c^T)+\mu$ is odd by (\ref{6.3}). And there holds
$b_{2j}=0$ for all $j\in\N_0$ by Lemma 2.6. Thus by (\ref{6.4}) we obtain
\be  B(d,h)(i(c^T) + p(c)) \ge -\sum_{2j-1=\mu-p(c)+1}^{i(c^T)+\mu}b_{2j-1}.
        \lb{6.5}\ee
Let $D=d(h+1)-2$. By Lemma 2.4 we have
$$  B(d,h) = -\frac{h(h+1)d}{2D}<0.  $$
Thus from (B) of Theorem 4.3, we have
\be    i(c^T)+\mu-(d-1)=i(c^T)+p(c)+dh-d-2 \in 2\N.  \lb{6.6}\ee
By (\ref{6.5}), (\ref{6.6}) and (\ref{2.12}) we obtain
\bea  i(c^T) + p(c)
&\le& -\frac{1}{B(d,h)}\sum_{2j-1=\mu-p(c)+1}^{i(c^T)+\mu}b_{2j-1}    \nn\\
&=& \frac{2D}{h(h+1)d}\left(\sum_{2j-1=1}^ {i(c^T)+\mu}b_{2j-1}-\sum_{2j-1=1}^{dh-2}b_{2j-1}\right).
             \lb{6.7}\eea
Note that here because $i(c^T)+p(c)\ge 2$ by (\ref{6.2}), we have
\be  i(c^T)+\mu = i(c^T)+p(c)+dh-3 \ge dh-1=d-1+(h-1)d. \lb{6.8}\ee

By Lemma 2.6 we have
\be \sum_{2j-1=1}^{i(c^T)+\mu}b_{2j-1} =
 \frac{h(h+1)d}{2D}\left(i(c^T)+\mu-(d-1)\right) - \frac{h(h-1)d}{4} + 1 +\ep_{d,h}(i(c^T)+\mu). \lb{6.9}\ee
On the other hand, because $dh-3< dh-1= d-1+(h-1)d$, by Lemma 2.6 we have
\bea \sum_{0\le 2j-1\le dh-3}b_{2j-1}
&=& \sum_{d-1\le 2j-1\le dh-3}\left([\frac{2j-1-(d-1)}{d}]+1\right)    \nn\\
&=& \sum_{d\le 2j\le dh-2}[\frac{2j}{d}]   \nn\\
&=& \sum_{\frac{d}{2}\le j\le \frac{dh}{2}-1}[\frac{j}{d/2}]   \nn\\
&=& \sum_{i=1}^{h-1}\sum_{j=\frac{id}{2}}^{\frac{(i+1)d}{2}-1}[\frac{j}{d/2}]   \nn\\
&=& \frac{d}{2}\sum_{i=1}^{h-1}i  \nn\\
&=& \frac{dh(h-1)}{4}.  \lb{6.10}\eea
Therefore we get
\bea
&& \sum_{0\le 2j-1\le i(c^T)+\mu}b_{2j-1} - \sum_{0\le 2j-1\le dh-3}b_{2j-1}   \nn\\
&&\qquad\quad = \frac{h(h+1)d}{2D}\left(i(c^T)+\mu-(d-1)\right) - \frac{h(h-1)d}{4}
             + 1 +\ep_{d,h}(i(c^T)+\mu) - \frac{dh(h-1)}{4} \nn\\
&&\qquad\quad = \frac{h(h+1)d}{2D}\left(i(c^T)+p(c)+dh-d-2\right) - \frac{dh(h-1)}{2}
             + 1 + \ep_{d,h}(i(c^T)+\mu). \quad \lb{6.11}\eea
Then (\ref{6.7}) becomes
$$  i(c^T) + p(c)
\le i(c^T)+p(c) + dh -d-2 + \frac{2D}{h(h+1)d}\left(1-\frac{dh(h-1)}{2}+\ep_{d,h}(i(c^T)+\mu)\right), $$
that is,
\bea  \ep_{d,h}(i(c^T)+\mu)
&\ge& \frac{h(h+1)d}{2D}\left(d+2 + \frac{(h-1)D}{h+1} - dh - \frac{2D}{h(h+1)d}\right)   \nn\\
&=& \frac{dh-(d-2)}{dh+(d-2)}.  \lb{6.12}\eea

Note that by (\ref{6.6}) we have
\be i(c^T)+\mu-(d-1) = i(c^T)+p(c)+dh-d-2 = i(c^T)+p(c)-2d + D. \lb{6.13}\ee
Let $\eta\in [0,D/2-1]$ be an integer such that
\be \frac{2\eta}{D} = \{\frac{i(c^T)+p(c)-2d}{D}\} = \{\frac{i(c^T)+\mu-(d-1)}{D}\}. \lb{6.14}\ee

By the definition (\ref{2.13}) of $\ep_{d,h}(i(c^T)+\mu)$ and (\ref{6.14}), we obtain
\bea  \ep_{d,h}(i(c^T)+\mu)
&=& \{\frac{D}{dh}\{\frac{i(c^T)+\mu-(d-1)}{D}\}\}
         - (\frac{2}{d}+\frac{d-2}{dh})\{\frac{i(c^T)+\mu-(d-1)}{D}\}   \nn\\
&&\qquad - h\{\frac{D}{2}\{\frac{i(c^T)+\mu-(d-1)}{D}\}\}
         - \{\frac{D}{d}\{\frac{i(c^T)+\mu-(d-1)}{D}\}\}    \nn\\
&=& \{\frac{2\eta}{dh}\} - (\frac{2}{d}+\frac{d-2}{dh})\frac{2\eta}{D}
      - h\{\frac{2\eta}{2}\} - \{\frac{2\eta}{d}\}   \nn\\
&=& \{\frac{2\eta}{dh}\} - (\frac{2}{d}+\frac{d-2}{dh})\frac{2\eta}{D} - \{\frac{2\eta}{d}\}  \nn\\
&\equiv& \ep(2\eta).   \lb{6.15}\eea

Now we claim
\be  \ep(2\eta) < \frac{dh-(d-2)}{dh+(d-2)}, \qquad \forall\; 2\eta\in [0,dh-2].  \lb{6.16}\ee

In fact, we write
\be  2\eta = pd + 2m \qquad \mbox{with some}\;\; p\in \N_0, \;\; 2m\in [0,d-2]. \lb{6.17}\ee
Then from $pd+2m=2\eta \le dh-2=(h-1)d+d-2$ we have
\be  p\in [0, h-1].  \lb{6.18}\ee
Therefore in this case we obtain
\bea \ep(2\eta)
&=& \frac{pd+2m}{dh} - (\frac{2}{d}+\frac{d-2}{dh})\frac{pd+2m}{D} - \frac{2m}{d}  \nn\\
&=& \frac{p}{h} - \frac{(2h+d-2)p}{hD} + \frac{2m}{dh} - \frac{(2h+d-2)2m}{dhD} - \frac{2m}{d}  \nn\\
&=& \frac{p}{h}(1-\frac{2h+d-2}{D}) + \frac{2m}{d}(\frac{1}{h} - \frac{2h+d-2}{hD} - 1)  \nn\\
&=& \frac{p(d-2) - 2mh}{D}  \nn\\
&\le& \frac{(h-1)(d-2)}{D}.  \lb{6.19}\eea
Now if (\ref{6.16}) does not hold, we then obtain
$$ \frac{dh-(d-2)}{D} \le \ep(2\eta) \le \frac{(h-1)(d-2)}{D}, $$
that is,
$$  dh - d +2 \le dh - d + 2 - 2h. $$
Because $h\ge 1$, this yields a contradiction and completes the proof of (\ref{6.16}).

If $d=2$, there holds $D-2=dh+d-4=dh-2$. Thus (\ref{6.16}) holds for any
integer $2\eta\in[0,D-2]$.

If $d\ge 4$, for any $2\eta\in [dh,D-2]$, write $2\eta = pdh + 2m$ for some $p\in\N_0$ and
$2m\in [0,dh-2]$. Then from $D-2=(h+1)d-4=hd+d-4$ we obtain $p\le 1$ and $2m\le d-4$.
Thus we have
\bea \ep(2\eta)
&=& \frac{2m}{dh} - (\frac{2}{d}+\frac{d-2}{dh})\frac{pdh+2m}{D} - \frac{2m}{d}  \nn\\
&=& \ep(2m) - (\frac{2}{d}+\frac{d-2}{dh})\frac{pdh}{D}  \nn\\
&\le& \ep(2m).  \lb{6.20}\eea

Therefore from (\ref{6.16}) and (\ref{6.20}), it yields
\be\ep(2\eta)<\frac{dh-(d-2)}{dh+(d-2)},\qquad\forall\,\eta\in [0,D/2-1].\lb{6.21}\ee
Together with (\ref{6.12}), (\ref{6.15}), and the
choice (\ref{6.14}) of $2\eta$, we then obtain
\be  \frac{dh-(d-2)}{dh+(d-2)} \le \ep_{d,h}(i(c^T)+\mu) = \ep(2\eta) < \frac{dh-(d-2)}{dh+(d-2)}.
            \lb{6.22}\ee
This contradiction completes the proof of Case 1.

\medskip

{\bf Case 2.} {\it $d\ge 2$ is odd and $h=1$.}

\medskip

In this case, $M$ is rationally homotopic to the sphere $S^d$. Note that $i(c^T)+\mu$
is even by (\ref{6.3}). Because $\mu-p(c)+1=d-2$, there holds $b_j=0$ for any
$j\le\mu-p(c)+1$ by Lemma 2.5. Thus by Theorem 5.5 and Lemma 2.5 we obtain
\bea  i(c^T) + p(c)
&\le& \frac{1}{B(d,1)}\sum_{2j=0}^{i(c^T)+\mu}b_{2j}        \nn\\
&\le& i(c^T) + p(c)  +d-3 - \frac{d-1}{2}\frac{2(d-1)}{d+1}        \nn\\
&=& i(c^T) + p(c) - \frac{4}{d+1}.   \lb{6.23}\eea
This is a contradiction.

\medskip

Now the proof of Theorem 1.1 is complete. \hfill\hb

\medskip

Now we give

{\bf Proof of Theorem 1.2.}  Here arguments are the same as in Section 7 of
\cite{LoD1}. For any reversible
Finsler as well as Riemannian metric $F$ on a compact manifold $M$, the energy
functional $E$ is symmetric on every loop $f\in \Lm M$ and its inverse curve
$f^{-1}$ defined by $f^{-1}(t)=f(1-t)$. Thus these two curves have the same
energy $E(f)=E(f^{-1})$ and play the same roles in the variational structure of
the energy functional $E$ on $\Lm M$. Specially, the $m$-th iterates $c^m$ and
$c^{-m}$ of a prime closed geodesic $c$ and its inverse curve $c^{-1}$ have
precisely the same Morse indices, nullities, and critical modules. Let
$n=n(c)=n(c^{-1})$. Then there holds
\be   \dim\ol{C}_*(E,c^m)=\dim\ol{C}_*(E,c^{-m}).    \lb{6.24}\ee
Thus if $c$ is the only geometrically distinct prime closed geodesic on $M$,
each entry in the diagram (\ref{5.54}) in the reversible case should be doubled.
So (\ref{5.46}) in Theorem 5.5 becomes
\be B(d,h)(i(c^T) + p(c)) + (-1)^{\mu+i(c^T)}2\ka
            = \sum_{j=\mu-p(c)+1}^{\mu+i(c^T)}(-1)^j b_j.  \lb{6.25}\ee
These changes bring no influence to our proofs in Section 6. Therefore our above
proof yields two geometrically distinct closed geodesics for reversible Finsler
metrics too. \hfill\hb

\bibliographystyle{abbrv}

\end{document}